\documentclass[hidelinks,onefignum,onetabnum]{siamart250211}

\usepackage{lipsum}
\usepackage{amsfonts}
\usepackage{graphicx}
\usepackage{epstopdf}

\usepackage{algorithmic}
\ifpdf
  \DeclareGraphicsExtensions{.eps,.pdf,.png,.jpg}
\else
  \DeclareGraphicsExtensions{.eps}
\fi

\newsiamremark{remark}{Remark}
\newsiamremark{hypothesis}{Hypothesis}
\crefname{hypothesis}{Hypothesis}{Hypotheses}
\newsiamthm{claim}{Claim}

\headers{Regularized multilevel methods}{N. Tsipinakis, P. Parpas}

\title{
Multilevel Regularized Newton Methods with Fast Convergence Rates
\thanks{Submitted to the editors July 2024.}
}

\author{Nick Tsipinakis\thanks{Department of Mathematics and Computer Science, UniDistance Suisse, Switzerland 
  (\email{nikolaos.tsipinakis@unidistance.ch}).}
\and Panos Parpas\thanks{Department of Computing, Imperial College, UK 
  (\email{panos.parpas@imperial.ac.uk}).}
}

\usepackage{amsopn}

\usepackage{mathtools}
\newcommand{\blue}{\textcolor{black}}
\newcommand{\red}{\textcolor{black}}

\usepackage{amssymb}

\usepackage{cases}

\ifpdf
\hypersetup{
  pdftitle={Multilevel Regularized Newton Methods with Fast Convergence Rates},
  pdfauthor={N. Tsipinakis, P. Parpas}
}
\fi

\newcommand{\R}{\mathbb{R}}
\newcommand{\N}{\mathbb{N}}
\newcommand{\w}{\mathbf}
\newcommand{\p}{\mathcal}
\newcommand{\q}{\mathrm}
\renewcommand{\d}{\mathsf{d}}
\newtheorem{assumption}{Assumption}
\usepackage{subcaption}

\begin{document}

\maketitle

\begin{abstract}
We introduce new multilevel methods for solving large-scale unconstrained optimization problems. Specifically, the philosophy of multilevel methods is applied to Newton-type methods that regularize the Newton sub-problem using second-order information from a coarse (low dimensional) sub-problem. The new \emph{regularized multilevel methods} provably converge from any initialization point and enjoy faster convergence rates than Gradient Descent. In particular, for arbitrary functions with Lipschitz continuous Hessians, we show that their convergence rate interpolates between the rate of Gradient Descent and that of the cubic Newton method. If, additionally, the objective function is assumed to be convex, then the proposed method converges with the fast $\mathcal{O}(k^{-2})$ rate. Hence, since the updates are generated using a \emph{coarse} model in low dimensions, the theoretical results of this paper significantly speed-up the convergence of Newton-type or preconditioned gradient methods in practical applications. Preliminary numerical results suggest that the proposed multilevel algorithms are significantly faster than current state-of-the-art methods. 
\end{abstract}

\begin{keywords}
Multilevel methods, Newton-type methods, unconstrained optimization
\end{keywords}

\begin{MSCcodes}
90C25, 90C26, 90C15
\end{MSCcodes}

\section{Introduction}
Let $f_H: \R^{N} \rightarrow \R$ be a twice differentiable function. In this work we are interested in solving the unconstrained optimization problem:
\begin{equation} \label{eq: optimization problem}
\w{x}_H^* = \underset{\w{x}_H \in \R^N}{\operatorname{arg \  min}} {f_H(\w{x}_H)}.
\end{equation}
The state-of-the-art second-order methods for solving (\ref{eq: optimization problem}) are the Newton method \cite{MR2142598}, Cubic Newton methods \cite{nesterov2008accelerating,nesterov2006cubic} and regularized Newton methods \cite{doikov2024gradient,mishchenko2023regularized}. Among them, the pure Newton method has a rich and powerful theory that shows that the pure Newton method has a quadratic convergence rate. Unfortunately, its rich theory can only be established locally, in a neighborhood of $\w{x}_H^*$. Imposing extra assumptions, such as strict convexity, it is possible to show global rates using a damping parameter (damped Newton method).
However, the damped Newton method has typically worse convergence rate than the $\p{O}(k^{-1})$ of Gradient Descent \cite{nesterov2018lectures}, where $k$ is the iteration number. The Cubic Newton method effectively addresses the theoretical limitations of Newton's method by adding a regularization term on the diagonal of the Hessian matrix. Assuming only the Lipschitz continuity of the Hessian matrix, the Cubic Newton method has the global state-of-the-art $\p{O} (\varepsilon^{-3/2})$ convergence rate (for gradient norm minimization), where $\varepsilon>0$ is the algorithm's tolerance. In addition, for general convex functions, the method enjoys the fast $\p{O}(k^{-2})$ rate (for function residual minimization). Both rates are better than those of Gradient Descent ($\p{O}(\varepsilon^{-2})$ and $\p{O}(k^{-1})$, respectively), justifying the extra computational cost imposed by forming the Hessian matrix and solving the resulting system of linear equations to obtain the Newton direction. Even though, these fast rates are achieved by solving a sub-problem at each iteration (which adds an extra computational cost to the Newton method), the Cubic Newton method is considered the method of choice when solving medium-sized non-convex problems. On the other hand, when minimizing general convex functions, it possible to show the same fast rates using gradient regularization \cite{doikov2024gradient,mishchenko2023regularized}. These methods do not need to solve the cubically regularized sub-problem and thus are faster than the Cubic Newton.

\red{    The main limitations of the classical second-order methods emerge when $N$ is very large, which typically arises in modern optimization problems. Under this regime, Newton-type methods either drastically slow-down or, even worse, they cannot be applied due to memory issues. Therefore, to take advantage of the valuable second-order information in large-scale optimization, multilevel or subspace methods offer an attractive alternative [1, 18, 19]. These methods effectively transfer the computations of the search directions in a \emph{coarse} level (subspace) where forming and inverting the Hessian matrix can be achieved efficiently. However, these methods, as we discuss in detail in the next section, suffer from the same limitations as the pure Newton method. In particular, they typically obtain worse global convergence rates than  Gradient Descent. It is even difficult to establish local convergence rates that are faster than $\mathcal{O}(k^{-1})$, and therefore their extra computational cost compared to Gradient Descent cannot be justified. \emph{The objective of this paper is to address the computational limitations of regularized Newton-type methods using a multilevel or randomized optimization framework. The proposed  approaches retain the robust theoretical guarantees of these methods, including rapid global convergence rates, even for non-convex problems, thereby maintaining their advantage over first-order methods.}}
    
\red{To achieve the desired rates, the deterministic multilevel approach requires computing a trial step at each iteration to dynamically assess efficiency each iterate according to a principled rule. Although this trial step may increase the computational cost relative to classical multilevel methods, the theory is important because it establishes convergence rates comparable to those of exact Newton methods and it provides a paradigm for designing reduced second-order methods with such fast convergence rates. Building on these insights, we develop a family of randomized regularized Newton methods that removes the need for a trial step while preserving the theoretical advantages of the deterministic multilevel method introducing a positive probability.}

\subsection{The Conventional Multilevel Method and other Related Work} \label{sec: conventional multilevel}
We present the main components of the conventional (Newton-type) multilevel method \cite{ho2019newton, tsipinakis2024multilevel,MR2587737}. In the multilevel literature the subscripts $H$ and $h$ are used to denote quantities in the fine (original) or coarse level, respectively, e.g., $f_H$ will be the fine model subject to minimization in problem (\ref{eq: optimization problem}) and $f_h: \R^{n} \rightarrow \R$, where $n < N$, will be the coarse model.
Further, we define $\w{P} \in \R^{N \times n}$ and $\w{R} \in \R^{n \times N}$ as the prolongation and restriction operators, respectively. These operators will be used to transfer the information between the coarse and fine levels. As such, given a point $\w{x}_{H,k}$, at iteration $k$, the initial point in the coarse level will be denoted as $\w{x}_{h,0} \coloneqq \w{R} \w{x}_{H,k}$. Then the coarse model is defined as follows:
\begin{equation} \label{eq: coarse model conventional multilevel}
    f_h(\w{x}_h) \coloneqq \langle \w{R}_k \nabla f_H (\w{x}_H), \w{x}_h - \w{x}_{h,0} \rangle + \frac{1}{2} \langle \w{R}_k \nabla^2 f_H (\w{x}_H) \w{P}_k (\w{x}_h - \w{x}_{h,0}), \w{x}_h - \w{x}_{h,0} \rangle.
\end{equation}
If the objective function is strongly convex, then $\w{R} \nabla^2 f (\w{x}_H) \w{P} \succ 0$, and thus the coarse model has a unique solution.
\blue{Let $\w{x}_{h}^*$ be its minimum point and define the vector $\w{d}_{h,k} \coloneqq \w{x}_{h}^* - \w{x}_{h,0}$. Note that  $\w{d}_{h,k}$ is a vector in $\R^n$ and needs to be \emph{prolongated} to the original space. In particular, the \emph{coarse direction} is defined by prolongating $\w{d}_{h,k}$, i.e., $  \w{d}_{H,k} \coloneqq \w{P}_k (\w{x}_{h}^* - \w{x}_{h,0})$}. Putting these concepts all together, the conventional multilevel method performs the following updates:
\begin{equation} \label{eq: updates conventional multilevel}
\begin{split}
    \w{x}_{H, k+1} & = \w{x}_{H, k} - t_{H,k} \w{P}_k \underset{\w{d}_h \in \R^n}{\operatorname{min}} \left\lbrace \langle \w{R}_k \nabla f (\w{x}_H), \w{d}_h \rangle + \frac{1}{2} \langle \w{R}_k \nabla^2 f (\w{x}_H) \w{P}\w{d}_h, \w{d}_h \rangle \right\rbrace  \\
    & =\w{x}_{H, k} - t_{H,k} \w{P}_k [\w{R}_k \nabla^2 f_H (\w{x}_{H,k}) \w{P}_k]^{-1} \w{R}_k \nabla f_H (\w{x}_{H,k}),
\end{split}
\end{equation}
where $t_{H,k} > 0$ is the step size parameter. \blue{The analysis of multilevel methods is typically performed for fixed $\w{R}_k$, say $\w{R}$, for all iterations $k$}.  It is known that the coarse directions may not always be effective. Ineffective directions occur when $\nabla f_H (\w{x}_H) \in \operatorname{null} (\w{R})$ which implies $\w{d}_{H} = 0$. This issue is addressed by alternating between the fine and coarse levels when computing the search directions according to some appropriate conditions. Let $\mu \in (0,\min(\|\w{R}\|, 1)]$ and $\varepsilon >0$. The coarse steps will be used if
\begin{equation} \label{eq: conditions}
            \| \w{R} \nabla f (\w{x}_{H,k}) \| \geq \mu  \| \nabla f (\w{x}_{H,k}) \| \quad \text{and}
        \quad \| \w{R} \nabla f (\w{x}_{H,k}) \| > \varepsilon.
\end{equation}
On the other hand, if the conditions above are violated, which imply $\w{d}_{H,k} \approx 0$, the conventional multilevel method computes the direction in the fine level, which is equivalent to performing the (damped) Newton method. In total, the cost of a single iteration for the multilevel method is $\p{O}(nN + n^3)$ which is much lower than the $\p{O}(N^2 + N^3)$ iteration cost of the Newton method.

Even though the conventional multilevel method is able to mitigate the computational limitations of second-order methods, by construction, it inherits all the theoretical drawbacks of the Newton method. That is, superior convergence rates to Gradient Descent can  only be established locally, while global rates are not better \cite{ho2019newton, tsipinakis2024multilevel}. Other works on multilevel methods fail to improve  the rate of convergence of Gradient Descent even locally \cite{gratton2008recursive,tsipinakis2023simba, MR2587737}. Another limitation, besides the unsatisfactory rates of convergence, is that
their updates (in \cref{eq: updates conventional multilevel}) are not well-defined for solving problem (\ref{eq: optimization problem}) in the absence of strong convexity. Even when convexity is present, multilevel or subspace methods based on randomization fail to improve the $\p{O}(k^{-1})$ convergence rate of Gradient Descent \cite{cartis2018global, cartis2022randomised}. \blue{Similar limitations are observed in several randomized sub-space methods. In particular, methods such as \cite{gower2019rsn,hanzely2023sketch} suffer from the same drawbacks as the conventional multilevel method since they require positive definite Hessians. Moreover, although the stochastic sub-space cubic Newton method \cite{hanzely2020stochastic} alleviates this assumption, requiring positive semi-definite Hessians, the $\p{O}(k^{-2})$ rate is attained only when $n=N$.}

\subsection{Contributions}

In view of the related work, we summarize the contributions of this paper:
\begin{enumerate}
    \item We introduce a family of regularized multilevel methods in order to speed-up the convergence of second-order methods for large-scale optimization.
    \item   The resulting regularized multilevel methods are analyzed for non-convex functions with Lipschitz continuous Hessians and for general convex functions. We show global convergence results in both cases. Furthermore, for non-convex functions with Lipschitz continuous Hessians, we show that the convergence rate of the proposed methods interpolates between the rate of Gradient Descent and the rate of the Cubic Newton method. In addition, for general convex functions, the regularized Newton-type multilevel method enjoys the fast $\p{O} (k^{-2})$ convergence rate of Cubic Newton method.
    \item When $\w{R}$ is a (particular) random matrix, we provide probabilistic results that allow the proposed methods to always perform iterations in the coarse level. Thus, such results further reduce the computational cost of the proposed methods. In the randomized case we show that it is possible to achieve the fast rate of \cite{mishchenko2023regularized}, under an appropriate choice of $\w{R}$ (\cref{def: P}) and specific problem structures. 
    \item We perform numerical experiments on real-world problems to show that the regularized multilevel methods achieve comparable rates compared to the Cubic Newton method. Moreover, we illustrate that the proposed methods are significantly faster in wall-clock time against both Gradient Descend and the Cubic Newton method.
\end{enumerate}

The rest of the paper is organized as follows: In \cref{sec: regularized multilevel methods} we describe the regularized multilevel methods and present the main algorithms. In \cref{sec: convergence analysis} we analyze the proposed methods and derive their global rates. In \cref{sec: implementation}, we discuss how to efficiently apply the proposed methods in practice and provide probabilistic convergence analysis. In \cref{sec: experiments} we present the numerical results.

\section{Regularized Multilevel Algorithms} \label{sec: regularized multilevel methods}

If \cref{eq: optimization problem} is not convex, then the coarse model in \cref{eq: coarse model conventional multilevel} of the conventional multilevel method is not suitable since the reduced Hessian matrix may not be positive definite or invertible. As a result, in order to obtain an effective coarse direction, one needs to replace the reduced Hessian matrix with an appropriate approximation. Let $\w{B}_h \coloneqq \w{B}_h(\w{x}_H) \in \R^{n \times n} $, be a symmetric positive semi-definite matrix and $\alpha_k > 0$. Then, we substitute $\w{R}_k \nabla^2 f_H (\w{x}_{H,k}) \w{P}_k$ in \cref{eq: coarse model conventional multilevel} with $\w{B}_{h, k} + \alpha_k \w{I}_n \succ 0$ (in \cref{sec: implementation} we will discuss possible candidates for $\w{B}_{h,k}$), and thus the update rule of the regularized multilevel method can be derived by minimizing the coarse model, which implies 
\begin{equation} \label{eq: update rule non convex coarse}
    \w{x}_{H, k+1} = \w{x}_{H, k} -  \w{P}_k (\w{B}_{h, k} + \alpha_k \w{I}_n)^{-1} \w{R}_k \nabla f_H (\w{x}_{H,k}).
\end{equation}
Note that, unlike the conventional method does, the regularized multilevel method does not require a step-size parameter. Instead, it adds a positive scalar to the diagonal of $\w{B}_k$. Therefore, a suitable $\alpha_k > 0$ should be added to ensure a sufficient decrease in the value of the objective function at each iteration. Since $\w{B}_{h, k}$ is an approximation of $\w{R}_k \nabla^2 f_H (\w{x}_{H,k}) \w{P}_k$ it is natural to define the error of this estimation
\begin{equation} \label{eq: definition s_h} 
    s_h \coloneqq \sup_{\w{x}_H} \| \w{B}_{h} - \w{R} \nabla^2 f_H (\w{x}_{H}) \w{P}\|.
\end{equation}

As with the standard multilevel method, one should alternate between the coarse and fine levels whenever $\w{R}_k \nabla f_H (\w{x}_{H,k}) = 0$ to ensure the progress of scheme (\ref{eq: update rule non convex coarse}). \blue{Note that, unlike \cref{eq: conditions}, where a fixed $\w{R}_k$ is used, here we allow $\w{R}_k$ to vary at each iteration, thereby enhancing the algorithm's applicability. To achieve fast convergence rates we introduce an additional condition alongside \cref{eq: conditions}. In particular, given a sequence of $\w{R}_k$, we use the following conditions to determine the space the search direction will be computed
\red{
\begin{equation} \label{eq: conditions new}
     \| \w{R}_k \nabla f (\w{x}_k) \| \geq \mu  \| \nabla f (\w{x}_k) \| \quad \text{and}
        \quad    \| \w{R}_k \nabla f (\w{x}_{k+1}) \| \geq \mu  \| \nabla f (\w{x}_{k+1}) \|,
\end{equation}
}
where $\mu \in (0, \min\{ \| \w{R}_{k-1} \|, \| \w{R}_k \|, 1\})$.
\red{Since the second inequality in \eqref{eq: conditions new} depends on $k+1$ iterate, evaluating it requires computing a trial step at each iteration. The new condition emerges from our theory and we will see later that it effectively ensures the algorithm's rapid convergence.}
If conditions \cref{eq: conditions new} are satisfied we use \cref{eq: update rule non convex coarse}, otherwise we use,}
\begin{equation} \label{eq: update rule non convex fine}
    \w{x}_{H, k+1} =  \w{x}_{H, k} -   (\w{B}_{H, k} + \alpha_k \w{I}_N)^{-1}  \nabla f_H (\w{x}_{H,k}),
\end{equation}
where $\w{B}_{H,k} \coloneqq \w{B}_{H,k} (\w{x}_H) \in \R^{N \times N}$ is a symmetric positive semi-definite approximation of the true Hessian matrix. \red{As a result, when the new condition is rejected, the method will compute both coarse and the full regularized Newton step at an iteration $k$. In such scenario the method forms expensive iterates. We address this issue in section 4.} 

As in \cref{eq: definition s_h}, we define the estimation error in the original space
\begin{equation} \label{eq: definition s_H}
    s_H \coloneqq \sup_{\w{x}_H} \| \w{B}_{H} -  \nabla^2 f_H (\w{x}_{H}) \|.
\end{equation}
On the other hand, if the problem is convex it is natural to employ the (reduced) Hessian matrix.
That is,
\begin{equation}\label{eq: update rule convex}
\w{x}_{H, k+1} =
\begin{cases} 
     \w{x}_{H, k} -  \w{P}_k (\w{R}_k \nabla^2 f_H (\w{x}_{H,k}) \w{P}_k + \alpha_k \w{I}_n)^{-1} \w{R}_k \nabla f_H (\w{x}_{H,k}),& \text{if \cref{eq: conditions new} holds} \\
      \w{x}_{H, k} -   (\nabla^2 f_H (\w{x}_{H,k}) + \alpha_k \w{I}_N)^{-1}  \nabla f_H (\w{x}_{H,k}), & \text{otherwise} 
\end{cases}
\end{equation}
\blue{Observe that if $n=N$ then we only use the second case in \cref{eq: update rule convex} above (since $\w{R}=\w{P}=\w{I}_N$).} Further, the reduced or exact Hessian matrix should be preferred in the convex setting because a faster convergence rate is expected due to the fact that $s_h = s_H = 0$. The regularized multilevel algorithms for non-convex and convex settings are presented in \cref{alg: non convex} and \cref{alg: convex}, respectively. Further, the algorithms require knowledge of either $s_H, s_h$ or the Lipschitz constant $L$ for computing $\alpha_k$, however these parameters are unknown in practice. Nevertheless, a simple line-search process can be applied. We will discuss these practical details in \cref{sec: implementation}.

\begin{algorithm}[]
\caption{Regularized multilevel algorithm - non-convex problems}
\begin{algorithmic}[1]
\STATE Input: $\w{x}_{0} \in \mathbb{R}^n, \w{R}_k \in \R^{n \times N}, \mu \in (0,1], L, s_h, s_H > 0$ \\ 
\FOR{$k = 1, \ldots$}
\IF{$    \red{ \| \w{R}_k \nabla f (\w{x}_k) \| \geq \mu  \| \nabla f (\w{x}_k) \| \ \text{and}
        \   \| \w{R}_k \nabla f (\w{x}_{k+1}) \| \geq \mu  \| \nabla f (\w{x}_{k+1}) \|}$}
\STATE Set $\alpha_k = s_h + \sqrt{\frac{L \|  \w{R}_k \nabla f (\w{x}_{H,k})\|}{2}}$
\STATE $\w{x}_{H, k+1} = \w{x}_{H, k} -  \w{P}_k (\w{B}_{h, k} + \alpha_k \w{I}_n)^{-1} \w{R}_k \nabla f_H (\w{x}_{H,k})$
\ELSE
\STATE Set $\alpha_k = s_H + \sqrt{\frac{L \|  \nabla f (\w{x}_{H,k})\|}{2}}$
\STATE $\w{x}_{H, k+1} =  \w{x}_{H, k} -   (\w{B}_{H, k} + \alpha_k \w{I}_n)^{-1}  \nabla f_H (\w{x}_{H,k}),$
\ENDIF
\ENDFOR
\end{algorithmic}
\label{alg: non convex}
\end{algorithm}

\begin{algorithm}[]
\caption{Regularized multilevel algorithm - convex problems}
\begin{algorithmic}[1]
\STATE Input: $\w{x}_{0} \in \mathbb{R}^n, \w{R}_k \in \R^{n \times N}, \mu \in (0,1], L > 0$ \\ 
\FOR{$k = 1, \ldots$}
\IF{$    \red{ \| \w{R}_k \nabla f (\w{x}_k) \| \geq \mu  \| \nabla f (\w{x}_k) \| \ \text{and}
        \   \| \w{R}_k \nabla f (\w{x}_{k+1}) \| \geq \mu  \| \nabla f (\w{x}_{k+1}) \|}$}
\STATE Set $\alpha_k = \sqrt{\frac{L \|  \w{R}_k \nabla f (\w{x}_{H,k})\|}{2}}$
\STATE $\w{x}_{H, k+1} = \w{x}_{H, k} -  \w{P}_k (\w{R}_k \nabla^2 f_H (\w{x}_{H,k}) \w{P}_k + \alpha_k \w{I}_n)^{-1} \w{R}_k \nabla f_H (\w{x}_{H,k})$
\ELSE
\STATE Set $\alpha_k = \sqrt{\frac{L \|  \nabla f (\w{x}_{H,k})\|}{2}}$
\STATE $\w{x}_{H, k+1} =  \w{x}_{H, k} -   (\nabla^2 f_H (\w{x}_{H,k}) + \alpha_k \w{I}_n)^{-1}  \nabla f_H (\w{x}_{H,k}),$
\ENDIF
\ENDFOR
\end{algorithmic}
\label{alg: convex}
\end{algorithm}

\section{Convergence Analysis} \label{sec: convergence analysis}
In this section we present the convergence analysis of Algorithms \ref{alg: non convex} and \ref{alg: convex} and derive the rates of convergence. Throughout this section, a vector $\w{x} \in \R^n$ will be denoted with bold lower-case letters while a matrix $\w{A} \in \R^{m \times n}$ with upper-case bold letters. Further, $\w{I}_p$, $p \in \N \setminus \{0\}$, denotes the $p \times p$ identity matrix. To simplify the notation, we discard the subscript $H$ from now onwards. In particular,  the quantities $\w{x}, f(\w{x}), \nabla f(\w{x})$ and $ \nabla^2 f(\w{x})$ will be assigned to $\w{x}_H, f_H(\w{x}_H), \nabla f_H(\w{x}_H)$ and $\nabla^2 f_H(\w{x}_H) $, respectively, unless noted otherwise. We use $\| \cdot\| $ to denote the standard Euclidean norm of a vector or a matrix. It holds $\| \w{x} \| = \sqrt{\langle \w{x}, \w{x} \rangle}$ and $\| \w{A}\| = \sqrt{\lambda_{\max(\w{A}^T \w{A})}}$ where $\lambda_{\max}$ corresponds to the largest eigenvalue of a square matrix \cite{MR2978290}. For two symmetric matrices $\w{B}, \w{D}$ we write $\w{B} - \w{D} \succeq 0$ if and only if $\langle\w{x} (\w{B} - \w{D}), \w{x} \rangle \geq 0.$
\subsection{Main inequalities and assumptions}
Our theory requires that the objective function has Lipschitz continuous Hessian matrices.
\begin{assumption} \label{ass: lipschitz cont}
    There exist a scalar $L>0$ such that for all $\w{x}, \w{y} \in \R^N$ it holds
    \begin{equation*}
        \| \nabla^2 f (\w{x}) - \nabla^2 f (\w{y}) \| \leq L \|\w{x} - \w{y} \|.
    \end{equation*}
\end{assumption}

We also make the following standard assumption on the operators.

\begin{assumption} \label{ass: P}
Let $k \in \N$.    It holds $\w{R}_k = \w{P}_k^T, \operatorname{rank} (\w{R}_k) = n$ and the sequence $(\|\w{R}_k\|)_{k \in \N}$ is bounded from above. 
\end{assumption}
The above assumption on $\w{P}$ and $\w{R}$ is typical and is satisfied by a number of choices that we discuss in \cref{sec: implementation}. Furthermore, \cref{ass: P} implies that there exist some $\omega \in (0, \infty)$ such that
    \begin{equation} \label{def: omega}
        \omega \coloneqq \sup \{ \| \w{R}_k\| \ | \ k \in \N \}.
    \end{equation}
We further define the following quantities:
\begin{align*}
    \hat{\lambda}(\w{x}_k) & \coloneqq [(\w{R}_k \nabla f(\w{x}_k))^T (\w{B}_{h,k} + \alpha_k \w{I}_n)^{-1} \w{R}_k \nabla f(\w{x}_k)]^{\frac{1}{2}} \\
    \w{d}_{h, k} &\coloneqq - (\w{B}_{h,k} + \alpha_k \w{I}_n)^{-1} \w{R}_k \nabla f(\w{x}_k)\\
    \w{d}_{H, k} &\coloneqq - \w{P}_k(\w{B}_{h,k} + \alpha_k \w{I}_n)^{-1} \w{R}_k \nabla f(\w{x}_k)
\end{align*}
Note that $\hat{\lambda}(\w{x}_k) \geq 0$ and it is a quantity analogous to the Newton decrement \cite{nesterov2018lectures}.
\begin{lemma} \label{lemma: identities}
    Consider the sequence $(\w{x}_{H,k})_{k \in \N}$ in \cref{eq: update rule non convex coarse}. Then 
    \begin{align}
        \langle \nabla f(\w{x}_k), \w{d}_{H, k} \rangle &= - \hat{\lambda}(\w{x}_k)^2  \label{eq: identity 1 app} \\
        \langle \nabla^2 f(\w{x}_k) \w{d}_{H, k}, \w{d}_{H, k} \rangle & \leq (s_h - \alpha_k) \|\w{d}_{h,k}\|^2 + \hat{\lambda}(\w{x}_k)^2 \label{eq: identity 2 app} 
    \end{align}
\end{lemma}

\begin{proof}
The equality in \cref{eq: identity 1 app} follows directly from the definitions of $\w{d}_{H,k}$ and $\hat{\lambda}(\w{x}_k)$. As for \cref{eq: identity 2 app}, we have that $\langle \nabla^2 f(\w{x}_k) \w{d}_{H, k}, \w{d}_{H, k} \rangle $
\begin{align*}
& = (\w{R}_k \nabla f(\w{x}_k))^T (\w{B}_{h,k} + \alpha_k \w{I}_n)^{-1} \w{R}_k \nabla^2 f(\w{x}_k) \w{P}_k (\w{B}_{h,k} + \alpha_k \w{I}_n)^{-1} \w{R}_k \nabla f(\w{x}_k) \\
& = (\w{R}_k \nabla f(\w{x}_k))^T (\w{B}_{h,k} + \alpha_k \w{I}_n)^{-1} (\w{R}_k \nabla^2 f(\w{x}_k) \w{P}_k - \w{B}_{h,k})(\w{B}_{h,k} + \alpha_k \w{I}_n)^{-1} \w{R}_k \nabla f(\w{x}_k)  \\
& \quad + (\w{R}_k \nabla f(\w{x}_k))^T (\w{B}_{h,k} + \alpha_k \w{I}_n)^{-1} \w{B}_{h,k} (\w{B}_{h,k} + \alpha_k \w{I}_n)^{-1} \w{R}_k \nabla f(\w{x}_k) \\
& = (\w{R}_k \nabla f(\w{x}_k))^T (\w{B}_{h,k} + \alpha_k \w{I}_n)^{-1} (\w{R}_k \nabla^2 f(\w{x}_k) \w{P}_k - \w{B}_{h,k})(\w{B}_{h,k} + \alpha_k \w{I}_n)^{-1} \w{R}_k \nabla f(\w{x}_k) \\
& \quad + (\w{R}_k \nabla f(\w{x}_k))^T (\w{B}_{h,k} + \alpha_k \w{I}_n)^{-1} (\w{B}_{h,k} + \alpha_k\w{I}_n) (\w{B}_{h,k} + \alpha_k \w{I}_n)^{-1} \w{R}_k \nabla f(\w{x}_k) \\ 
& \quad - \alpha_k \| \w{d}_{h,k}\|^2 \\
&\leq s_h \| \w{d}_{h,k}\|^2 + \hat{\lambda}(\w{x}_k)^2 - \alpha_k \| \w{d}_{h,k}\|^2
\end{align*}
where for the last inequality we used $\langle \w{A} \w{x}, \w{x} \rangle \leq \|\w{A}\| \|\w{x}\|^2$, in particular,
\begin{align*}
    &(\w{R}_k \nabla f(\w{x}_k))^T (\w{B}_{h,k} + \alpha_k \w{I}_n)^{-1} (\w{R}_k \nabla^2 f(\w{x}_k) \w{P}_k - \w{B}_{h,k})(\w{B}_{h,k} + \alpha_k \w{I}_n)^{-1} \w{R}_k \nabla f(\w{x}_k) \\ & \leq  s_h \| \w{d}_{h,k}\|^2,
\end{align*}
and thus (\ref{eq: identity 2 app}) has been proved.
\end{proof}

The following bounds will be useful in our analysis several times:
\begin{equation} \label{ineq: upper bound d_H}
    \begin{split}
    \| \w{d}_{h,k} \| \leq \| (\w{B}_{h,k} + \alpha_k \w{I}_n)^{-1}\| \| \w{R}_k \nabla f(\w{x}_k)\| \leq 
    \frac{\| \w{R}_k \nabla f(\w{x}_k)\|}{\alpha_k}, \ \| \w{d}_{H,k} \| \leq \omega \| \w{d}_{h,k} \|, 
    \end{split}
\end{equation}
where for the result on the left we have used that if $0 \prec \w{A} \preceq \w{B}$, then $\|\w{A}\| \leq \| \w{B} \|$. The next result justifies the choice of $\alpha_k$ in \cref{alg: non convex}.
\begin{lemma} \label{lemma: lower bound d_h}
    Let Assumptions \ref{ass: lipschitz cont} and \ref{ass: P} hold and suppose that the sequence $(\w{x}_k)_{k \in \N}$ is generated by \cref{eq: update rule non convex coarse}. If $\alpha_k \geq s_h + \sqrt{\frac{\omega^3 L \| \w{R}_k \nabla f(\w{x}_k)\|}{2}}$, then
    \begin{equation*}
        \| \w{R}_k \nabla f(\w{x}_{k+1})\| \leq 2 \alpha_k \|\w{d}_{h,k} \|.
    \end{equation*}
\end{lemma}

\begin{proof}
    It holds
    \begin{align*}
    \| \w{R}_k \nabla f(\w{x}_{k+1})\|
        & = \| \w{R}_k \nabla f(\w{x}_{k+1}) - \w{R}_k \nabla f(\w{x}_{k}) - (\w{B}_{h,k} + \alpha_k \w{I}_n) \w{d}_{h,k}\| \\
        & \leq \| \w{R}_k \nabla f(\w{x}_{k+1}) - \w{R}_k \nabla f(\w{x}_{k}) - \w{B}_{h,k} \w{d}_{h,k}\| + \alpha_k \| \w{d}_{h,k}\| \\
        & \leq\| \w{R}_k \nabla f(\w{x}_{k+1}) - \w{R}_k \nabla f(\w{x}_{k}) - \w{R}_k \nabla^2 f(\w{x}_{k}) \w{P}_k \w{d}_{h,k}\| \\
        & \qquad \qquad + \| (\w{R}_k \nabla^2 f(\w{x}_{k}) \w{P}_k - \w{B_k}) \w{d}_{h,k}\| + \alpha_k \| \w{d}_{h,k}\| \\
        & \leq \omega \| \nabla f(\w{x}_{k} + \w{d}_{H,k}) - \nabla f(\w{x}_{k}) - \nabla^2 f(\w{x}_{k}) \w{d}_{H,k}\| + (s_h + \alpha_k) \|\w{d}_{h,k} \| \\
        & = \omega \left\| \int_0^1 \left( \nabla^2 f(\w{x}_{k} + t \w{d}_{H,k}) - \nabla^2 f(\w{x}_{k}) \right) \d t \  \w{d}_{H,k} \right\| + (s_h + \alpha_k) \|\w{d}_{h,k} \| \\
        & \leq  \omega \int_0^1 \| \nabla^2 f(\w{x}_{k} + t \w{d}_{H,k}) - \nabla^2 f(\w{x}_{k}) \| \d t  \ \|  \w{d}_{H,k}\| + (s_h + \alpha_k) \|\w{d}_{h,k} \| \\
        & \overset{\text{Ass. }\ref{ass: lipschitz cont}}{ \leq} \frac{\omega L \|\w{d}_{H,k} \|^2}{2} + (s_h + \alpha_k) \|\w{d}_{h,k} \| \\
        & \overset{(\ref{ineq: upper bound d_H})}{\leq}  \frac{\omega^3 L \|\w{d}_{h,k} \|^2}{2} + (s_h + \alpha_k) \|\w{d}_{h,k} \| \\
        & \overset{(\ref{ineq: upper bound d_H})}{\leq} \left(\frac{\omega^3 L \| \w{R}_k \nabla f(\w{x}_{k})\| }{2 \alpha_k} + s_h + \alpha_k \right)\|\w{d}_{h,k} \| \leq 2 \alpha_k \|\w{d}_{h,k} \|. \\
    \end{align*}
\end{proof}
 The next result shows that the sequence $(f(\w{x}_k))_{k \in \N}$, where $\w{x}_k$ as in \cref{eq: update rule non convex coarse}, is non-increasing given $\alpha_k$ in \cref{alg: non convex}.
\begin{lemma} \label{lemma: f reduction lambda}
    Let Assumptions \ref{ass: lipschitz cont} and \ref{ass: P} hold and suppose that the sequence $(\w{x}_k)_{k \in N}$ is generated by \cref{eq: update rule non convex coarse}. If
    $\alpha_k \geq s_h + \sqrt{\frac{\omega^3 L \| \w{R}_k \nabla f(\w{x}_k)\|}{2}}$, then
    \begin{equation*}
        f(\w{x}_{k+1}) \leq f(\w{x}_k )- \frac{\hat{\lambda}(\w{x}_k)^2}{2}.
    \end{equation*}
\end{lemma}

\begin{proof}
    By \cref{ass: lipschitz cont} it can be shown that (see \cite{nesterov2018lectures} for a proof)
        \begin{align*}
        f(\w{x}_{k+1}) 
        & \leq f(\w{x}_{k}) + \langle  \nabla f(\w{x}_{k}), \w{x}_{k+1} - \w{x}_{k}\rangle + \frac{1}{2} \langle  \nabla^2 f(\w{x}_{k}) (\w{x}_{k+1} - \w{x}_{k}), \w{x}_{k+1} - \w{x}_{k} \rangle \\ & \qquad \quad \ \ + \frac{L}{6} \|\w{x}_{k+1} - \w{x}_{k}\|^3 \\
        & = f(\w{x}_{k}) + \langle  \nabla f(\w{x}_{k}), \w{d}_{H,k}\rangle + \frac{1}{2} \langle  \nabla^2 f(\w{x}_{k}) \w{d}_{H,k}, \w{d}_{H,k} \rangle + \frac{L}{6} \|\w{d}_{H,k}\|^3.
        \end{align*}
        Using \cref{eq: identity 1 app} and \cref{eq: identity 2 app} we obtain
        \begin{align*}
            f(\w{x}_{k+1}) - f(\w{x}_{k})
            & \leq - \hat{\lambda}(\w{x}_k)^2 + \frac{1}{2}[ (s_h - \alpha_k) \|\w{d}_{h,k}\|^2 + \hat{\lambda}(\w{x}_k)^2] + \frac{L}{6} \|\w{d}_{H,k}\|^3 \\
            & \overset{(\ref{ineq: upper bound d_H})}{\leq} - \frac{\hat{\lambda}(\w{x}_k)^2}{2} +\frac{1}{2} (s_h - \alpha_k) \|\w{d}_{h,k}\|^2 + \frac{\omega^3 L}{6} \|\w{d}_{h,k}\|^3 \\
            & \overset{(\ref{ineq: upper bound d_H})}{\leq} - \frac{\hat{\lambda}(\w{x}_k)^2}{2} +\frac{1}{2} (s_h - \alpha_k) \|\w{d}_{h,k}\|^2 + \frac{\omega^3 L \|\w{R}_k \nabla f(\w{x}_{k})\|}{6 \alpha_k} \|\w{d}_{h,k}\|^2 \\
            & = -\frac{\hat{\lambda}(\w{x}_k)^2}{2} - \frac{1}{2} \left( \alpha_k - s_h - \frac{\omega^3 L \|\w{R}_k \nabla f(\w{x}_{k})\|}{3 \alpha_k} \right) \|\w{d}_{h,k}\|^2\\
            & \leq -\frac{\hat{\lambda}(\w{x}_k)^2}{2},
        \end{align*}
        where the last inequality follows from our assumption on $\alpha_k$ which ensures that $\alpha_k - s_h - \frac{\omega^3 L \|\w{R}_k \nabla f(\w{x}_{k})\|}{3 \alpha_k}>0$.
\end{proof}
Moreover, from $    \w{B}_{h,k} + \alpha_k \w{I}_n \succeq \alpha_k \w{I}_n,$ we have that
\begin{align*}
         (\w{B}_{h,k} + \alpha_k \w{I}_n)^{-1} \preceq \frac{1}{\alpha_k} \w{I}_n
        & \iff   \w{x}^T \left((\w{B}_{h,k} + \alpha_k \w{I}_n)^{-1} - \frac{1}{\alpha_k} \w{I}_n \right) \w{x}  \leq 0, \   \text{for all} \ \w{x} \in \R^n \\ &\iff \w{x}^T (\w{B}_{h,k} + \alpha_k \w{I}_n)^{-1}  \w{x}  \leq \frac{\|\w{x}\|^2}{\alpha_k},
        \  \text{for all} \ \w{x} \in \R^n         
\end{align*}
Thus, setting $\w{x}\coloneqq (\w{B}_{h,k} + \alpha_k \w{I}_n)^{-\frac{1}{2}} \w{R}_k \nabla 
f(\w{x}_{k})$ we take that
\begin{equation} \label{ineq: lower bound lamda}
    \|\w{d}_{h,k}\|^2 \leq \frac{\hat{\lambda}(\w{x}_k)^2}{\alpha_k}
\end{equation}
The results of this section will be used below to derive the rates of convergence for \cref{alg: non convex} and \cref{alg: convex}.

\subsection{The Non-Convex Case} \label{sec: non-convex analysis}
In this section we show global convergence of a subsequence of \cref{alg: non convex} to a stationary point of $f$ and derive its convergence rate. We will need the following definition 
\begin{equation*}
    (\q{g}_k)_{k \in \N} \coloneqq (\| \nabla f(\w{x}_{k}) \|)_{k \in \N}.
\end{equation*}
Next, we combine the results in (\ref{ineq: lower bound lamda}) and \cref{lemma: f reduction lambda} to obtain a reduction in the value of the objective function in terms of the norm of the gradient.
\begin{corollary} \label{cor: f reduction norm g}
    Let Assumptions \ref{ass: lipschitz cont} and \ref{ass: P} hold and suppose that the sequence $(\w{x}_k)_{k \in N}$ is generated by \cref{eq: update rule non convex coarse}. If
    $\alpha_k \geq s_h + \sqrt{\frac{\omega^3 L \| \w{R}_k \nabla f(\w{x}_k)\|}{2}}$, then
    \begin{equation*}
        f(\w{x}_{k+1}) \leq f(\w{x}_k )- \frac{\mu^2 \| \nabla f(\w{x}_{k+1})\|^2}{8 \alpha_k}.
    \end{equation*}    
\end{corollary}

\begin{proof}
By \cref{lemma: f reduction lambda} and \cref{lemma: lower bound d_h} we have that 
\begin{align*}
f(\w{x}_k ) - f(\w{x}_{k+1}) \geq  \frac{\hat{\lambda}(\w{x}_k)^2}{2} \overset{(\ref{ineq: lower bound lamda})}{\geq}  \frac{\alpha_k \|\w{d}_{h,k}\|^2}{2} \geq \frac{\|\w{R}_k \nabla f(\w{x}_{k+1})\|^2}{8\alpha_k} \geq \frac{\mu^2 \| \nabla f(\w{x}_{k+1})\|^2}{8\alpha_k}
\end{align*}
where the last inequality holds because we generate the sequence $(\w{x}_k)_{k \in \N}$ \red{in the coarse
level and thus \eqref{eq: conditions new} is satisfied.} 
\end{proof}

\begin{remark}
If we additionally assumed that $(\q{g}_k)_{k \in \N}$ is bounded, then the result of \cref{cor: f reduction norm g} could be used to prove that $(\q{g}_k)_{k \in \N}$ is a null sequence. That is, by telescoping the above bound and using $\alpha_k \coloneqq s_h + \sqrt{\omega^3 \frac{L}{2}} \q{g}_i^{\frac{1}{2}}$ we get 
\begin{equation*}
    f(\w{x}_0) - f(\w{x}^*) \geq f(\w{x}_0) - f(\w{x}_{k}) \geq \sum_{i=0}^{k-1} \frac{\mu^2 \q{g}_{i+1}^2}{8\left( s_h + \sqrt{\omega^3 \frac{L}{2}} \q{g}_i^{\frac{1}{2}} \right) }.
\end{equation*}
Therefore, as $k \rightarrow \infty$, the series converges, which together with the boundness $(\q{g}_k)_{k \in \N}$ of imply that $\lim_{k \rightarrow \infty} \q{g}_k = 0$. Thus, the sequence $(\w{x}_k)_{k \in \N}$, constructed in \cref{eq: update rule non convex coarse}, converges to a stationary point of problem (\ref{eq: optimization problem}) from any initialization point $\w{x}_0$.
\end{remark}

Unfortunately, even under this assumption, the convergence rate of $(\q{g}_k)_{k \in \N}$ is unknown. Nevertheless, we can derive a rate of convergence for one of its subsequences, namely
\begin{equation*}
    \q{g}_{k}^* \coloneqq \min_{0 \leq i \leq k} \| \nabla f(\w{x}_{i}) \|.
\end{equation*}
The proof of the next result is parallel to that of \cite{doikov2024gradient}. 
\begin{lemma} \label{lemma: core non-convex}
    Let $k \in \N$ and $\varepsilon>0$ such that $\q{g}_i \geq \varepsilon$  for all $i \in \{0,1, \ldots, k\}$. Then 
\begin{equation*}
    k \leq \left\lfloor \frac{8(f(\w{x}_0) - f^*)}{\mu^2} \left(\sqrt{\frac{L}{2}} \frac{\omega^2}{\varepsilon^{\frac{3}{2}}} + \frac{s_h}{\varepsilon^2}
        \right) + 2\ln{\frac{\q{g}_0}{\varepsilon}} \right\rfloor.
\end{equation*}
\end{lemma}

\begin{proof}
    From \cref{cor: f reduction norm g} and since $\| \w{R}_k \nabla f(\w{x}_{k}) \| \leq \omega \q{g}_k$ we take
    \begin{align*}
        f(\w{x}_k) - f(\w{x}_{k+1}) 
        & \geq \frac{\mu^2 \q{g}_{k+1}^2}{8\left(s_h + \omega^2 \sqrt{\frac{L}{2}} \q{g}_k^\frac{1}{2}\right)} \\
        & =  \left( \frac{\q{g}_{k+1}}{\q{g}_{k}} \right)^2 \frac{\mu^2}{8 \left(s_h \q{g}_{k}^{-2} + \omega^2 \sqrt{\frac{L}{2}} \q{g}_{k}^{-\frac{3}{2}}\right)} \\
        & \geq  \left( \frac{\q{g}_{k+1}}{\q{g}_{k}} \right)^2 \frac{\mu^2}{8 \left(s_h \varepsilon^{-2} + \omega^2 \sqrt{\frac{L}{2}} \varepsilon^{-\frac{3}{2}}\right)},
    \end{align*}
    where in the last inequality we used the assumption $\q{g}_k \geq \varepsilon$. Define
    \begin{equation*}
        C_\varepsilon \coloneqq \frac{\mu^2}{8 \left(s_h \varepsilon^{-2} + \omega^2 \sqrt{\frac{L}{2}} \varepsilon^{-\frac{3}{2}}\right)}.
    \end{equation*}
    Telescoping the previous bound we obtain
    \begin{align*}
        f(\w{x}_0) - f(\w{x}^*) \geq f(\w{x}_0) - f(\w{x}_{k}) \geq C_\varepsilon \sum_{i=0}^{k-1} \left( \frac{\q{g}_{i+1}}{\q{g}_{i}}\right)^2.
    \end{align*}
    Using the inequalities
    \begin{equation*}
        \exp(x) \geq 1 + x \quad \text{and} \quad \frac{1}{k}\sum_{i=1}^{k} x_i \geq \left(\prod_{i=1}^k x_i \right)^{\frac{1}{k}}
    \end{equation*}
    we have that
    \begin{align*}
        f(\w{x}_0) - f(\w{x}^*) 
        & \geq C_\varepsilon \sum_{i=0}^{k-1} \left( \frac{\q{g}_{i+1}}{\q{g}_{i}}\right)^2 \geq C_\varepsilon k \prod_{i=1}^{k-1} \left( \frac{\q{g}_{i+1}}{\q{g}_{i}}\right)^\frac{2}{k} = C_\varepsilon k \left( \frac{\q{g}_{k}}{\q{g}_{0}}\right)^\frac{2}{k} \geq 
        C_\varepsilon k \left( \frac{\varepsilon}{\q{g}_{0}}\right)^\frac{2}{k} \\
        & = C_\varepsilon k \exp \left(\frac{2}{k} \ln \frac{\varepsilon}{\q{g}_{0}}\right) \geq C_\varepsilon k \left(1 + \frac{2}{k} \ln \frac{\varepsilon}{\q{g}_{0}}\right)
    \end{align*}
    and thus, solving for $k$, the result follows.
\end{proof}

The result in \cref{lemma: core non-convex} is key for obtaining the convergence rate. Effectively, it reveals the first iteration $k$, (smallest $k \in \N$), such that $\q{g}_k < \varepsilon$. It is now easy to obtain the convergence rate of the subsequence $(\q{g}_k^*)_{k \in \N}$.
\begin{theorem} \label{thm: coarse non-convex}
    Let Assumptions \ref{ass: lipschitz cont} and \ref{ass: P} hold. Suppose also that the sequence $(\w{x}_k)_{k \in \N}$ is generated by (\ref{eq: update rule non convex coarse}) and 
    \begin{equation*}
        \alpha_k = s_h + \sqrt{\frac{\omega^3 L\| \w{R}_k \nabla f(\w{x}_{k}) \|}{2}},
    \end{equation*}
    where $s_h$ is defined in \cref{eq: definition s_h}.
    Then for every $\varepsilon > 0$ and for
    \begin{equation*}
    K_\varepsilon \coloneqq \left\lfloor \frac{8(f(\w{x}_0) - f^*)}{\mu^2} \left(\sqrt{\frac{L}{2}} \frac{\omega^2}{\varepsilon^{\frac{3}{2}}} + \frac{s_h}{\varepsilon^2}
        \right) + 2\ln{\frac{\q{g}_0}{\varepsilon}} \right\rfloor + 1
\end{equation*}
    we have that $\q{g}_k^* < \varepsilon$ for all $k \geq K_\varepsilon$.
\end{theorem}

\begin{proof}
    Let $\varepsilon>0$. From \cref{lemma: core non-convex} it holds directly that if $k \geq K_\varepsilon$, then there exists at least one $i \in \{0,1, \ldots, k\}$ such that $\q{g}_i < \varepsilon$. Thus, by the definition of $\q{g}_k^*$ we have that $\q{g}_{K_\varepsilon}^* < \varepsilon$. Since also $(\q{g}_k^*)_{k \in \N}$ is non-increasing, it holds $\q{g}_k^* < \varepsilon$ for all $k > K_\varepsilon$, which completes the proof.
\end{proof}

\cref{thm: coarse non-convex} shows that if all the iterations were performed in the coarse level, then we could reach tolerance $\varepsilon > 0$ after at least $K_\varepsilon$ iterations. The convergence rate in \cref{thm: coarse non-convex} effectively interpolates between the convergence rate of gradient descent and cubic Newton algorithms (plus a logarithmic constant) according to the values of $s_h$ and $\mu$. For instance, if $s_h$ is large then the method converges in $\mathcal{O}(\varepsilon^{-2})$ iterations, while if $s_h=0$ the method converges in $\mathcal{O}(\varepsilon^{-\frac{3}{2}})$ iterations. 
As for $\mu$, recall that it is a user-defined parameter which controls the number of steps taken between the coarse and fine levels. Therefore, the user has access to the rate of convergence by tuning $\mu$, nevertheless  that large values in $\mu$ are more likely to require (expensive) iterations with the full model. We complete the convergence analysis of \cref{alg: non convex} by providing a similar result, but now all iterations are performed in the fine level.

\begin{theorem} \label{thm: fine non-convex}
        Let \cref{ass: lipschitz cont} hold. Suppose also that the sequence $(\w{x}_k)_{k \in \N}$ is generated by \cref{eq: update rule non convex fine} and 
    \begin{equation*}
        \alpha_k = s_H + \sqrt{\frac{L\| \nabla f(\w{x}_{k}) \|}{2}},
    \end{equation*}
    where $s_H$ is defined in \cref{eq: definition s_H}. Then for every $\varepsilon > 0$ and 
    \begin{equation*}
    K_\varepsilon \coloneqq \left\lfloor 8(f(\w{x}_0) - f^*) \left(\sqrt{\frac{L}{2}} \frac{1}{\varepsilon^{\frac{3}{2}}} + \frac{s_H}{\varepsilon^2}
        \right) + 2\ln{\frac{\q{g}_0}{\varepsilon}} \right\rfloor + 1
\end{equation*}
    we have that $\q{g}_k^* < \varepsilon$ for all $k \geq K_\varepsilon$.
\end{theorem}

The proof of the theorem can be obtained by simply replacing $s_h$ with $s_H$ and $\w{R}$ with $\w{I}_N$ in \cref{thm: coarse non-convex}, and thus it is omitted. Hence, Theorems \ref{thm: coarse non-convex} and \ref{thm: fine non-convex} provide the complete picture of the convergence rate of \cref{alg: non convex}. Further, we emphasize that this rate is better than that of Gradient Descent, which requires $\mathcal{O}(\frac{M}{\varepsilon^2})$ iterations, where $M>0$ such that $\| \nabla^2 f(\w{x}) \| \leq M$, for all $\w{x} \in \R^N$. This is because $s_h$ and $s_H$ can become arbitrarily small while $M$ is typically large. 
 \subsection{The Convex Case}

In this section we additionally assume that problem \cref{eq: optimization problem} is convex and analyze \cref{alg: convex} with the new assumption. Recall that the results in \cref{lemma: identities}, \cref{lemma: lower bound d_h} and \cref{lemma: f reduction lambda} require only Assumptions \ref{ass: lipschitz cont} and \ref{ass: P}. Since here the reduced Hessian is positive semi-definite the results of the aforementioned lemmas hold for $\w{B}_{h,k} \coloneqq \w{R}_k \nabla^2 f(\w{x}_{k}) \w{P}_k$. In this case it holds $s_h = 0$. We immediately obtain the following results.

\begin{corollary} \label{cor: f reduction norm g convex}
    Let $f$ be a convex function and  Assumptions \ref{ass: lipschitz cont} and \ref{ass: P} hold. Suppose also that the sequence $(\w{x}_k)_{k \in N}$ is generated by the first equation in \cref{eq: update rule convex}. If
    $\alpha_k \geq \sqrt{\frac{\omega^3 L \| \w{R}_k \nabla f(\w{x}_k)\|}{2}}$, then the following inequality are satisfied
        \begin{align*}
        & \| \w{R}_k \nabla f(\w{x}_{k+1})\| \leq 2 \alpha_k \|\w{d}_{h,k} \|, \\
        & f(\w{x}_{k+1}) \leq f(\w{x}_k )- \frac{1}{2} \hat{\lambda}(\w{x}_k)^2,\\
        & f(\w{x}_{k+1}) \leq f(\w{x}_k )- \frac{\mu^2 \| \nabla f(\w{x}_{k+1})\|^2}{8 \alpha_k}.
        \end{align*}
\end{corollary}
Let $f^* \coloneqq f(\w{x}^*)$ and further define
    \begin{equation} \label{def: R}
        R \coloneqq \sup_{\w{x} \in \R^N} \{ \|\w{x} - \w{x}^* \| : f(\w{x}) \leq f(\w{x}_0)\}.
    \end{equation}
Applying the proof technique of \cite[Theorem 1]{mishchenko2023regularized} to our results we show $\mathcal{O} (\frac{1}{k^2})$ convergence rate of \cref{alg: convex}. 

\begin{theorem} \label{thm: rate convex coarse}
    Let $f$ be a convex function and assumptions \ref{ass: lipschitz cont} and \ref{ass: P} hold. Suppose also that $R<\infty$, the sequence $(\w{x}_k)_{k \in \N}$ is generated by \cref{alg: convex} and set
    \begin{equation*}
        \alpha_k = \sqrt{\frac{\omega^3 L \| \w{R}_k \nabla f(\w{x}_k)\|}{2}}.
    \end{equation*}
    Then, there exists $K_0 \in \N$ such that convergence rate of \cref{alg: convex} is given by
    \begin{equation}
        f(\w{x}_k) - f^* \leq \frac{128^2 \omega^4 R^3 L}{2 \mu^4 \left(1 + \frac{3k}{4(\log_{2} \frac{\omega}{\mu} + 3)} \right)^2}, \ \ \text{for all} \ k \geq K_0\label{ineq: rate thm convex}
    \end{equation}
\end{theorem}

\begin{proof}
    Let us fix an arbitrary $k \in \N$. 
    First, we analyze the case where \cref{alg: convex} generates the updates always in the subspace, i.e.,  $(\w{x}_k)_{k \in \N}$ is generated by \cref{eq: update rule non convex coarse}. Then, by \cref{cor: f reduction norm g convex} and the fact that $\| \w{R}_k \nabla f(\w{x}_k)\| \leq \omega \| \nabla f(\w{x}_k)\|$ we have
    \begin{equation} \label{ineq: final f reduction in theorem} 
        f(\w{x}_{k}) - f(\w{x}_{k+1} ) \geq \frac{\mu^2 \| \nabla f(\w{x}_{k+1})\|^2}{8 \| \w{R}_k \nabla f(\w{x}_k)\|^{\frac{1}{2}} \sqrt{\frac{\omega^3 L }{2}}} \geq   \frac{\mu^2 \q{g}_{k+1}^2}{8 \omega^2  \sqrt{\frac{ L }{2}}\q{g}_{k}^{\frac{1}{2}}} ,
    \end{equation}
    where $\q{g}_{k} \coloneqq \| \nabla f(\w{x}_{k})\|$. Since the sequence $(f(\w{x}_m))_{m \in \N}$ is non-increasing, then by \cref{def: R} it holds $\|\w{x}_m - \w{x}^* \| \leq R$, for all $m \in \N$. Thus, by convexity and the Cauchy Schwartz inequality we have that
    \begin{equation} \label{ineq: convexity}
        f(\w{x}_k) - f^* \leq \langle \nabla f(\w{x}_k), \w{x}_k - \w{x}^* \rangle \leq R \q{g}_k.
    \end{equation}
    As in the proof of \cite[Theorem 1]{mishchenko2023regularized}, we define the following index sets: $\mathcal{I}_{\infty} \coloneqq \{ m \in \N: \q{g}_{m+1} \geq \frac{1}{4} \q{g}_m \}$ and $\mathcal{I}_{k} \coloneqq \{ m \in \mathcal{I}_\infty: m \leq k \}$ and the rate of convergence depends on the size of $\mathcal{I}_{k}$. 
    
    Let $D_k \coloneqq \frac{3k}{2(\log_{2}\frac{\omega}{\mu} + 3)}$. By the definition of $\mu$ in (\ref{eq: conditions}) we get $0< D_k \leq \frac{k}{2}$, and thus $D_{k+1} \geq D_k$ which implies $\lim_{k \rightarrow \infty } D_k = \infty$.
    
    \textbf{Case 1:} $|\mathcal{I}_k| \geq D_k$. Then, $\mathcal{I}_\infty$ contains infinitely many elements. Thus there exists a subsequence  $(m_\ell)_{\ell \in \N}$ such that $ \mathcal{I}_\infty = \{m_0, m_1, \ldots \}$.  
    Therefore, from (\ref{ineq: final f reduction in theorem}) we get 
    \begin{equation*}
        f(\w{x}_{m_{\ell}}) - f(\w{x}_{m_{\ell}+1} ) \geq \frac{\mu^2 \q{g}_{m_{\ell}+1}^2}{8 \omega^2  \sqrt{\frac{ L }{2}}\q{g}_{m_\ell}^{\frac{1}{2}}} \geq \frac{\mu^2 }{128 \omega^2  \sqrt{\frac{ L }{2}}} \q{g}_{m_{\ell}}^\frac{3}{2}.
    \end{equation*}
    Further, $m_{\ell +1} \geq m_{\ell} + 1 $ and since $(f(\w{x}_m))_{m \in \N}$ is non-increasing we have that
        \begin{equation*}
        f(\w{x}_{m_{\ell}}) - f(\w{x}_{m_{\ell+1}} )  \geq \frac{\mu^2 \q{g}_{m_{\ell}+1}^2}{8 \omega^2  \sqrt{\frac{ L }{2}}\q{g}_{m_\ell}^{\frac{1}{2}}} \geq \frac{\mu^2 }{128 \omega^2  \sqrt{\frac{ L }{2}}} \q{g}_{m_{\ell}}^\frac{3}{2}.
    \end{equation*}
    Next, let $C \coloneqq \frac{\mu^2 }{128 \omega^2 R^\frac{3}{2}  \sqrt{\frac{ L }{2}}}$ and combine the above inequality with (\ref{ineq: convexity}) to obtain
    \begin{align*}
        & f(\w{x}_{m_{\ell+1}}) - f(\w{x}_{m_{\ell}})  \leq - C [f(\w{x}_{m_{\ell}}) - f^*]^\frac{3}{2} \iff \\
        & \frac{f(\w{x}_{m_{\ell+1}}) - f^*}{C} \leq \frac{f(\w{x}_{m_{\ell}}) - f^*}{C} -  (f(\w{x}_{m_{\ell}}) - f^*)^{\frac{1}{2}} (f(\w{x}_{m_{\ell+1}} - f^*)        
    \end{align*}
    Then, using the notation $\Delta_\ell \coloneqq \sqrt{\frac{f(\w{x}_{m_{\ell}}) - f^*}{C}}$ we take $\Delta_{\ell+1}^2 \leq \Delta_{\ell}^2 - C^\frac{3}{2}\Delta_{\ell}^3$ and note that $\Delta_\ell \leq \frac{1}{C^\frac{3}{2}}$, for all $\ell \in \N$. Using this we obtain the following inequality
    \begin{equation*}
        \frac{1}{\Delta_{\ell +1}} - \frac{1}{\Delta_{\ell}} = \frac{\Delta_{\ell} - \Delta_{\ell+1}}{\Delta_{\ell}\Delta_{\ell+1}}= 
        \frac{\Delta_{\ell}^2 - \Delta_{\ell+1}^2}{\Delta_{\ell+1}\Delta_{\ell}(\Delta_{\ell}+\Delta_{\ell+1})} \geq C^{\frac{3}{2}}\frac{\Delta_{\ell}^3}{2\Delta_{\ell}^3} = \frac{C^\frac{3}{2}}{2}.
    \end{equation*}
    Next, by telescoping the last inequality we have that 
    \begin{equation*}
        \frac{1}{\Delta_{\ell}} \geq \frac{1}{\Delta_0} +  \ell \frac{C^{\frac{3}{2}}}{2} \geq C^\frac{3}{2} + D_k \frac{C^{\frac{3}{2}}}{2}, 
    \end{equation*}
    where the last inequality holds from the assumption $|\mathcal{I}_k| \geq D_k$ and the lower bound on $\Delta_0$. Inverting this inequality we can  derive the convergence rate for the first case. First we take
    \begin{equation*}
        f(\w{x}_{m_{\ell}}) - f^* \leq \frac{1}{[(1 + \frac{D_k}{2})C]^2}, \ \  \text{for all} \  \ell \in \N.
    \end{equation*}
    Then, since $|\mathcal{I}_k| \geq D_k$, there exist $\tilde{\ell} \in \N$ such that $ m_{\tilde{\ell}} \in \mathcal{I}_k$. Thus,  $m_{\tilde{\ell}} \leq k$ and $f(\w{x}_k) \leq f(\w{x}_{m_{\tilde{\ell}}})$. Putting all these together we have that
        \begin{equation}
        f(\w{x}_k) - f^* \leq f(\w{x}_{m_{\tilde{\ell}}}) - f^* \leq \frac{128^2 \omega^4 R^3 L}{2 \mu^4 \left(1 + \frac{3k}{4(\log_{2} \frac{\omega}{\mu} + 3)} \right)^2}.
    \end{equation}

    \textbf{Case 2:} $|\mathcal{I}_k| \coloneqq M < D_k$. For each $m \in \N$, we have two possibilities for $\q{g}_m$: \textbf{(a)}: If $m \notin \mathcal{I}_k \cap \{1,2, \ldots, k\}$, then $m \notin \mathcal{I}_\infty$ and thus $\q{g}_{m} < \frac{1}{4}\q{g}_{m-1}$. If, on the other hand, $m \in \mathcal{I}_k$, then we always have
    $\q{g}_{m} \leq  \frac{\omega}{\mu} 2 \q{g}_{m-1}$, since from \cref{cor: f reduction norm g convex} it holds that
    \begin{equation*}
        \mu \q{g}_{m+1} \overset{(\ref{eq: conditions})}{\leq} \| \w{R}_m \nabla f(\w{x}_{m+1})\| \leq 2 \alpha_m \|\w{d}_{h,m} \| \overset{(\ref{ineq: upper bound d_H})}{\leq} 2 \|\w{R}_m \nabla f(\w{x}_{m}) \| 
        \leq 2 \omega \q{g}_{m}.
    \end{equation*}
    Then, combining the two inequalities for $\q{g}_{m}$ with convexity (\ref{ineq: convexity}) we have that
    \begin{equation*}
       f(\w{x}_k) - f^* \leq R \q{g}_k \leq  \frac{(\frac{2\omega}{\mu})^M }{4^{k-M}} R \q{g}_0 \leq \frac{R}{2^{\frac{k}{2}}}\q{g}_0,
    \end{equation*}
    where the last inequality holds since $M < D_k$, which completes the analysis for our cases.

    Note that, in the second case, \cref{alg: convex} (using coarse steps) achieves an exponential convergence rate. Hence, since this is a much faster convergence, there exists $K_0 \in \N$ such that
    \begin{equation*}
        \frac{R}{2^{\frac{k}{2}}}\q{g}_0 \leq \frac{128^2 \omega^4 R^3 L}{2 \mu^4 \left(1 + \frac{3k}{4(\log_{2} \frac{\omega}{\mu} + 3)} \right)^2}, \ \ \text{for all} \ k \geq K_0. 
    \end{equation*}
    This proves \cref{ineq: rate thm convex}. Further, if $(\w{x}_k)_{k \in \N}$ is generated by the second equation in \cref{eq: update rule convex}, then the rate follows from \cite[Theorem 1]{mishchenko2023regularized}. This is faster than \cref{ineq: rate thm convex} for all $k$, and hence, we conclude that the overall convergence rate of \cref{alg: convex} is given by \cref{ineq: rate thm convex}.
    \end{proof}

\cref{thm: rate convex coarse} shows that \cref{alg: convex} achieves the fast $\mathcal{O}(k^{-2})$. In particular, in the first $K_0$ steps the algorithm may converge exponentially before it enters the $\mathcal{O}(k^{-2})$ rate. In comparison to \cref{thm: coarse non-convex}, the multilevel algorithm, in the convex setting, requires $\mathcal{O}(\varepsilon^{-1/2})$ iterations to reach a tolerance $\varepsilon>0$ which is a significant improvement to the best possible rate for non-convex functions. Furthermore, as in \cref{thm: coarse non-convex}, the user-defined parameter $\mu$ is proportional to the convergence rate of the algorithm. The smaller the parameter $\mu$ is selected by the user (in order to allow for more coarse steps), the slower the rate in \cref{thm: rate convex coarse} becomes. Hence, there is a trade-off between the fast rate and the number of coarse steps. It is possible to combine very fast rates with large number of coarse steps if additionally one assumes specific problem structures, for more details see \cite{tsipinakis2024multilevel, tsipinakis2023multilevel}.

\section{Line Search and Analysis with Random Operators} \label{sec: implementation}

In this section we focus on efficient implementation of Algorithms \ref{alg: non convex} and \ref{alg: convex}. Regarding \cref{alg: non convex} and the choice of $\w{B}_{h,k}$, we consider the following approximations for the reduced Hessian:
\begin{enumerate}
    \item If $n$ is selected small, then the eigenvalue decomposition of $\w{R}_k \nabla^2 f(\w{x}_k) \w{P}_k$ can be applied. Let $\w{U}_n \in \R^{n \times n}$ and $\Lambda \in \R^{n \times n}$ such that 
    $ \w{R}_k \nabla^2 f(\w{x}_k) \w{P}_k = \w{U}_n \Lambda_n \w{U}_n^T$, where $\w{U}_n$ is the matrix of eigenvectors with orthonormal columns and $\Lambda_n \coloneqq \operatorname{diag} (\lambda_1(\w{x}), \lambda_2(\w{x}), \ldots, \lambda_n(\w{x}))$ is a diagonal matrix consisting of the associated eigenvalues, and we assume that $\lambda_1(\w{x}) \geq \lambda_2(\w{x}) \geq \ldots \geq \lambda_n(\w{x})$. Then we define
    \begin{equation}\label{eq:scen1}
        \w{B}_{h} \coloneqq \w{U}_n |\Lambda_n| \w{U}_n^T, \quad  |\Lambda_n |\coloneqq \operatorname{diag} (|\lambda_1(\w{x})|, |\lambda_2(\w{x})|, \ldots, |\lambda_n(\w{x})|)  
    \end{equation}
    \label{scenario 1}
    \item Alternatively, to reduce the computational cost of approach \ref{scenario 1}, one may compute only the first $r < n$ eigenvectors and eigenvalues. This can be efficiently achieved using the power method described in \cite{halko2011finding} \footnote{Python implementation at: \url{https://pytorch.org/docs/stable/generated/torch.svd_lowrank.html}}. Then, the matrix $\w{B}_h$ is formed as follows.
    Let $\w{U}_r \in \R^{n \times r}$ and $ \Lambda_r \in \R^{r \times r}$ be the results of the power method in \cite{halko2011finding}. As with $\Lambda_n$, we define $|\Lambda_r|$ by taking absolute values and 
        \begin{equation} \label{eq:scen2}
        \w{B}_{h} \coloneqq \w{U}_r |\Lambda_r| \w{U}_r^T.
    \end{equation}
    Then, in practice, we compute $(\w{B}_{h,k} + \alpha_k\w{I}_n)^{-1}$ using the Woodbury identity
    \begin{equation*}
        (\w{B}_{h,k} + \alpha_k \w{I}_n)^{-1} = \frac{1}{\alpha_k} \w{I}_n - \frac{1}{\alpha_k^2} \w{U}_{r,k} \left(|\Lambda_r| + \frac{1}{\alpha_k} \w{U}^T_{r,k} \w{U}_{r,k} \right)^{-1} \w{U}^T_{r,k}
    \end{equation*}
        \label{scenario 2}
    \item In the last scenario, we regularize the reduced Hessian matrix by adding the absolute value of the minimum eigenvalue to its diagonal to obtain a positive semi-definite matrix:
    \begin{equation}\label{eq:scen3}
        \w{B}_{h,k} \coloneqq \w{R}_k \nabla^2 f(\w{x}_k) \w{P}_k + |\lambda_{\min}(\w{R}_k \nabla^2 f(\w{x}_k) \w{P}_k)| \w{I}_n,
    \end{equation}
    To achieve this, one can either compute the full eigenvalue decomposition and proceed as in \eqref{eq:scen1} or estimate only the smallest eigenvalue which is typically faster. \footnote{Python implementation at: \url{https://pytorch.org/docs/stable/generated/torch.lobpcg.html}}
        \label{scenario 3}
\end{enumerate}

Using the various definitions of $\w{B}_{h,k}$, we examine cases and reveal problem structures where the rate of \cref{alg: non convex} is superior to Gradient Descent. First, recall the rate shown in \cref{thm: coarse non-convex} and the rate of Gradient Descent, i.e., $\p{O}(\frac{s_h}{\varepsilon^{2}})$ and $\p{O}(\frac{M}{\varepsilon^{2}})$, respectively, where $M>0$ is such that $-M \w{I}_N \preceq \nabla^2 f(\w{x}_k) \preceq M \w{I}_N$, for all $\w{x} \in \R^N$. Then, for the algorithm implemented with  \eqref{eq:scen3}, we have
$s_h  \geq |\lambda_{\min}(\w{R}_k \nabla^2 f(\w{x}_k) \w{P}_k)|$ but $M \geq \max \{|\lambda_{\min}(\w{R}_k \nabla^2 f(\w{x}_k) \w{P}_k)|, |\lambda_{\max}(\w{R}_k \nabla^2 f(\w{x}_k) \w{P}_k)| \}$. Therefore, if problem \cref{eq: optimization problem} is ill-conditioned, then $M > s_h$ and typically $M \gg s_h$. Similarly, one can show that the algorithm implemented using \eqref{eq:scen2} is suitable for problems where there is gap between the $r$ and the $r+1$ eigenvalues, see also \cite{tsipinakis2023multilevel}.

Moreover, recall that algorithms \ref{alg: non convex} and \ref{alg: convex} require knowledge of parameters $L, S_h$ or $S_H$. Since these parameters are often unknown in practice, to estimate them, we apply a simple line search procedure similar to that used for Cubic Newton \cite{nesterov2018lectures}. We first describe the line search strategy for \cref{alg: convex} which requires only an estimation of $L$. The process requires an initial estimation $L_0>0$. Subsequently, to find $L_k$, the current value is repeatedly increased by a factor of two until the inequality in \cref{lemma: f reduction lambda} is satisfied. Then $L_{k+1}$ is initialized by reducing $L_k$ by a factor of two, and $L_{k+1}$ can be at least as low as $L_0$. The process is given in \cref{alg: inner loop convex}. We present the steps in the coarse level only for simplicity. 
\begin{algorithm}
\caption{\cref{alg: convex} with line search}
\begin{algorithmic}[1]
\STATE Input: $\w{x}_0, \ L_0 > 0,  \ \w{R}_k,  \ \mu \in (0, \min(\|\w{R}_k\|, 1)), \  \varepsilon \geq 0$ \\ 
\FOR{$k = 0,1, \ldots$}
\STATE $\hat{\lambda}(\w{x}_k)^2 = - \langle \nabla f (\w{x}_k), \w{d}_{H,k} \rangle$, where $\w{d}_{H,k}$ as in (\ref{eq: update rule convex})
\STATE $\w{x}_\text{step} = \w{x}_{k} -  \w{P}_k (\w{R}_k \nabla^2 f (\w{x}_{k}) \w{P}_k + \alpha_k \w{I}_n)^{-1} \w{R}_k \nabla f (\w{x}_{k})$
\STATE $i_k = 1$
\WHILE{$f (\w{x}_\text{step}) > f (\w{x}_k) - \frac{\hat{\lambda}(\w{x}_k)^2}{2}$} \label{step: exit condition}
\STATE $\alpha_k =  \sqrt{\frac{2^{i_k} L_k \|  \w{R}_k \nabla f (\w{x}_{H,k})\|}{2}}$
\STATE $\w{x}_\text{step} = \w{x}_{k} -  \w{P}_k (\w{R}_k \nabla^2 f (\w{x}_{k}) \w{P}_k + \alpha_k \w{I}_n)^{-1} \w{R}_k \nabla f (\w{x}_{k})$
\STATE $i_k = i_k + 1$
\ENDWHILE
\STATE $\w{x}_{k+1} = \w{x}_\text{step}$
\STATE $L_{k+1} = \max \{L_0, 2^{i_k - 1} L_k \}$ \label{step: inner loops}
\ENDFOR
\end{algorithmic}
\label{alg: inner loop convex}
\end{algorithm}

We can derive an upper bound on the number of inner loops for \cref{alg: inner loop convex}. 
From step \ref{step: inner loops} of \cref{alg: inner loop convex} we have that $i_k = \log_{2} \frac{L_{k+1}}{L_k} +1$, and thus assuming that the algorithms runs for a total number of $T$ iterations, we obtain
\begin{equation*}
    \sum_0^{T} i_k = T + 1 + \log_{2} \frac{L_{T+1}}{L_0} \leq T + 1 + \log_{2} \frac{L}{L_0},
\end{equation*}
where in the last inequality we have used that $L_k \leq L $ for all $k\in \N$. As for the non-convex variant, along with $L_k$, we introduce $s_k>0$ for estimating $s_h$. We will be updating these parameters simultaneously using the same principle, see \cref{alg: inner loop non convex}. Therefore, similar to \cref{alg: inner loop convex},  the total number of inner loops for the non-convex case given by
\begin{equation*}
    \sum_0^{T} i_k  \leq T + 1 + \max \{ \log_{2} \frac{L}{L_0}, \log_{2} \frac{s_h}{s_0}\}.
\end{equation*}

Moreover, several strategies have been proposed for selecting the operators $\w{R}_k$ and $\w{P}_k$. Among the most popular are the cyclic or randomized rules because convergence can be attained without ever applying the checking process (\ref{eq: conditions}) \cite{cartis2022randomised, drineas2005nystrom, gittens2011spectral,tsipinakis2024multilevel}. For the latter rule, at each iteration, the entries of operator $\w{R}_k$ are independently sampled from a distribution $\mathcal{D}$ such that \cref{ass: P} is satisfied. Of particular interest are distributions that impose sparsity on $\w{R}_k$ since matrix multiplications can be achieved more efficiently. Specifically, forming $\w{R}_k$ by sampling rows of the identity matrix significantly reduces the complexity of the iterates.

\begin{definition} \label{def: P}
    Let the $S_N \coloneqq \{1,2, \ldots, N \}$. Sample $n < N$ elements from $S_N$ without replacement from uniform distribution and denote $S_n \coloneqq \{s_1, s_2, \ldots, s_n \}$. Then, the $i^{\text{th}}$ row of $\w{R}_k$, $i \in \{1,2, \ldots, n \}$, is the $s_i$ row of the identity matrix $\w{I}_N$ and $\w{P}_k \coloneqq \w{R}_k$.  
\end{definition}

Using \cref{def: P} to form $\w{R}_k$, we no longer need to perform matrix-vector multiplications since $\w{R}_k$ now acts as a sampling operator. For instance, the product $\w{R} \w{y}$ for some fixed $\w{y} \in \R^N$ is equivalent to sampling $n$ out of $N$ elements from $\w{y}$. Similarly, one can efficiently compute the derivative of $\w{R}_k \nabla f (\w{x}_k)$ to obtain the reduced Hessian matrix (for more details see \cite{tsipinakis2024multilevel}). In a similar fashion, one can form the rows of $\w{R}_k$ from a selection of the rows of $\w{I}_N$ in a cyclic manner (deterministically). For instance, we can find all combinations of $N$ over $n$, which are $q=\binom{N}{n}$ in number, to cyclically update $\w{R}_1, \w{R}_2, \ldots, \w{R}_q$. Additionally, beyond the efficient implementation, forming $\w{R}_k$ through $\w{I}_N$ yields some nice properties for the iterates of the proposed algorithms that we aim to exploit in order to improve the current convergence rates.

\begin{algorithm}[t]
\caption{\cref{alg: non convex} with line search}
\begin{algorithmic}[1]
\STATE Input: $L_0, s_0 > 0$ \\ 
\FOR{$k = 0,1, \ldots$}
\STATE $\hat{\lambda}(\w{x}_k)^2 = - \langle \nabla f (\w{x}_k), \w{d}_{H,k} \rangle$, where $\w{d}_{H,k}$ as in (\ref{eq: update rule non convex coarse})
\STATE $\w{x}_\text{step} = \w{x}_{k} -  \w{P}_k (\w{R}_k \nabla^2 f (\w{x}_{k}) \w{P}_k + \alpha_k \w{I}_n)^{-1} \w{R}_k \nabla f (\w{x}_{k})$
\STATE $i_k = 1$
\WHILE{$f (\w{x}_{k+1}) > f (\w{x}_k) - \frac{\hat{\lambda}(\w{x}_k)^2}{2}$}
\STATE Set $\alpha_k = 2^{i_k} s_k + \sqrt{\frac{2^{i_k} L_k \|  \w{R}_k \nabla f (\w{x}_{H,k})\|}{2}}$
\STATE $\w{x}_\text{step} = \w{x}_{k} -  \w{P}_k (\w{B}_{k} + \alpha_k \w{I}_n)^{-1} \w{R}_k \nabla f (\w{x}_{k})$
\STATE $i_k = i_k + 1$
\ENDWHILE
\STATE $\w{x}_{k+1} = \w{x}_\text{step}$
\STATE Set $L_{k+1} = \max \{L_0, 2^{i_k - 1} L_k \}$ and $s_{k+1} = \max \{s_0, 2^{i_k - 1} s_k \}$
\ENDFOR
\end{algorithmic}
\label{alg: inner loop non convex}
\end{algorithm}

\subsection{Convergence Analysis with Random Operators}

In this section we present the convergence analysis of \cref{alg: non convex} and \cref{alg: convex} when $\w{R}_k$ is constructed according to \cref{def: P}. Hence, due to randomness, we are no longer required to perform the expensive fine step. Additionally, the convergence rates will now hold with some probability. We use the notation of \cref{sec: convergence analysis} and we further denote
\begin{equation*}
    r_k : = \|\w{x}_{k+1} -  \w{x}_{k}\|.
\end{equation*}

Moreover, by \cref{def: P}, we directly obtain the following results:
\begin{equation} \label{eq: core results for R}
    \w{R}_k \w{P}_k = \w{I}_n \quad \text{and} \quad \|\w{R}_k\| = 1, \ \text{for all} \  k \in \N.
\end{equation}
From this we also get $\omega=1$, where $\omega$ is defined in \cref{def: omega}. Given the properties of $\w{R}_k$ above we can prove the following identities.

\begin{lemma}
    Let $k \in \N$. For $\w{B}_{h,k} \succeq 0$ the following equalities are satisfied 
    \begin{align}
        -& \langle (\w{B}_{h,k} + \alpha_k \w{I}_n)\w{R}_k \w{d}_{H,k}, \w{R}_k \w{d}_{H,k} \rangle = \langle \w{R}_k \nabla f(\w{x}_k), \w{R}_k \w{d}_{H,k}\rangle \label{eq: identity 1}  \\ 
        & \langle  \nabla f(\w{x}_{k}), \w{d}_{H,k} \rangle =  \langle   \w{R}_k \nabla f(\w{x}_{k}), \w{R}_k \w{d}_{H,k} \rangle \label{eq: identity 2} \\
        & \langle  \nabla^2 f(\w{x}_{k})\w{d}_{H,k}, \w{d}_{H,k} \rangle =  \langle  \w{R}_k \nabla^2 f(\w{x}_{k})\w{P}_k \w{R}_k\w{d}_{H,k}, \w{R}_k \w{d}_{H,k} \rangle    \label{eq: identity 3}\\
        & r_k = \| \w{R}_k (\w{x}_{k+1} - \w{x}_{k})\| = \| \w{x}_{k+1} - \w{x}_{k}\|, \label{eq: identity 4}
    \end{align}
    where $\w{d}_{H,k} = \w{x}_{k+1} - \w{x}_{k}$.
\end{lemma}

\begin{proof}
    We prove the above identities for $\w{x}_{k+1} - \w{x}_{k}$ in place of $\w{d}_{H,k}$.
    Multiplying (\ref{eq: update rule non convex coarse}) with $\w{R}_k$  we have that
    \begin{equation}
        \w{R}_k(\w{x}_{k+1} - \w{x}_{k}) = -(\w{B}_{h,k} + \alpha_k \w{I}_n)^{-1} \w{R}_k \nabla f(\w{x}_k) \label{eq: identity 5}
    \end{equation}
    and thus, 
    \begin{align*}
        - & (\w{B}_{h,k} + \alpha_k \w{I}_n)\w{R}_k (\w{x}_{k+1} - \w{x}_{k})  = \w{R}_k \nabla f(\w{x}_k) \iff \\
          & \langle (\w{B}_{h,k} + \alpha_k \w{I}_n)\w{R}_k (\w{x}_{k+1} - \w{x}_{k}), \w{R}_k (\w{x}_{k+1} - \w{x}_{k}) \rangle = \langle \w{R}_k \nabla f(\w{x}), \w{R}_k (\w{x}_{k+1} - \w{x}_{k})\rangle,
    \end{align*}
    as claimed. For \cref{eq: identity 2} we take
    \begin{align*}
       \langle \nabla f(\w{x}_{k}), \w{x}_{k+1} - \w{x}_{k} \rangle & = -\langle \nabla f(\w{x}_{k}), \w{P}_k (\w{B}_{h,k} + \alpha_k \w{I}_n)^{-1} \w{R}_k \nabla f(\w{x}) \rangle \\
       & = - \langle \w{R}_k \nabla f(\w{x}_{k}),  (\w{B}_{h,k} + \alpha_k \w{I}_n)^{-1} \w{R}_k \nabla f(\w{x}) \rangle \\
       & \stackrel{(\ref{eq: identity 5})}{=}  \langle \w{R}_k \nabla f(\w{x}_{k}), \w{R}_k (\w{x}_{k+1} - \w{x}_{k}) \rangle
    \end{align*}
    Similar arguments can be used to prove \cref{eq: identity 3}.
    Finally, for \cref{eq: identity 4} we get
    \begin{align*}
        \| \w{R}_k (\w{x}_{k+1} - \w{x}_{k})\| & \leq \| \w{x}_{k+1} - \w{x}_{k}\| = \|\w{P}_k (\w{B}_{h,k} + \alpha_k \w{I}_n)^{-1} \w{R}_k \nabla f(\w{x}_k)\| \\ & \leq 
        \|(\w{B}_{h,k} + \alpha_k \w{I}_n)^{-1} \w{R}_k \nabla f(\w{x})\| \stackrel{(\ref{eq: identity 5})}{=} \| \w{R}_k (\w{x}_{k+1} - \w{x}_{k})\|
    \end{align*}
    which completes the proof.
\end{proof}

Furthermore, we will be using the following upper bound for $r_k$.
\begin{equation} \label{ineq: upper bound r_k}
\begin{split}
    r_k \leq \| (\w{B}_{h,k} + \alpha_k \w{I}_n)^{-1}\| \| \w{R}_k \nabla f(\w{x}_k)\| \leq \frac{\| \w{R}_k \nabla f(\w{x}_k)\|}{\alpha_k}.
\end{split}
\end{equation}

\subsubsection{The Non-Convex Case} \label{sec: non-convex probabilities}
Here we show probabilistic convergence rate of \cref{alg: non convex} with \cref{def: P}. Next, we derive a new bound on $\alpha_k$.

\begin{lemma} \label{lemma: lower bound rk}
    Let \cref{ass: lipschitz cont} and $\w{R}_k$ as in \cref{def: P}. Suppose also that the sequence $(\w{x}_k)_{k \in \N}$ is generated by \cref{eq: update rule non convex coarse}. If $\alpha_k \geq s_h + \sqrt{\frac{L \| \w{R}_k \nabla f(\w{x}_k)\|}{2}}$, then
    \begin{equation*}
        \| \w{R}_k \nabla f(\w{x}_{k+1})\| \leq 2 \alpha_k r_k.
    \end{equation*}
\end{lemma}

\begin{proof}
Let us denote $\w{d}_k\coloneqq - \w{P} (\w{B}_{h,k} + \alpha_k \w{I}_n)^{-1} \w{R}_k \nabla f(\w{x}_k)$. Then from \cref{eq: identity 5} we take that
\begin{equation} \label{eq: identity proof lower bound rk}
    \w{P}_k \w{R}_k (\w{x}_{k+1} - \w{x}_{k}) = \w{d}_k
\end{equation}
The proof follows that of \cref{lemma: lower bound d_h} and thus we report only the key steps.
It holds $\| \w{R}_k \nabla f(\w{x}_{k+1}) \| \overset{(\ref{eq: identity 5})}{=} $
\begin{align*}
    & \qquad \  \| \w{R}_k \nabla f(\w{x}_{k+1}) - \w{R}_k \nabla f(\w{x}_{k}) -  (\w{B}_{h,k} + \alpha_k \w{I}_n) \w{R}_k (\w{x}_{k+1} - \w{x}_{k})\| \\
    &\overset{\cref{eq: core results for R}}{\leq} \| \nabla f(\w{x}_{k+1}) - \nabla f(\w{x}_{k}) - \nabla^2 f(\w{x}_{k}) \w{P}_k \w{R}_k (\w{x}_{k+1} - \w{x}_{k}) \| + (s_h + \alpha_k) r_k \\
    & \overset{(\ref{eq: identity proof lower bound rk})}{=} \| \nabla f(\w{x}_{k} + \w{d}_k) - \nabla f(\w{x}_{k}) - \nabla^2 f(\w{x}_{k}) \w{d}_k \| + (s_h + \alpha_k) r_k \\
    & \leq \int_0^1 \| \nabla^2 f(\w{x}_{k} + t\w{d}_k) - \nabla^2 f(\w{x}_{k})\| \d t \ r_k + (s_h + \alpha_k) r_k \\
    & \overset{\text{Ass. }\ref{ass: lipschitz cont}}{\leq} \frac{L r_k^2}{2} + (s_h + \alpha_k) r_k \\
    &  \overset{(\ref{ineq: upper bound r_k})}{\leq} \left(\frac{L \| \w{R}_k \nabla f(\w{x}_{k})\|}{2\alpha_k} +  s_h + \alpha_k \right)r_k\\
    &\leq 2 \alpha_k r_k
\end{align*}
where the last inequality holds from the assumption on $\alpha_k$.  
\end{proof}


\begin{lemma} \label{lemma: f reduction definition R}
    Let Assumptions \ref{ass: lipschitz cont} and $\w{R}_k$ as in \cref{def: P}. Suppose also that the sequence $(\w{x}_k)_{k \in \N}$ is generated by \cref{eq: update rule non convex coarse}. Then for any $\alpha_k \geq s_h + \sqrt{\frac{L \| \w{R}_k \nabla f(\w{x}_k)\|}{2}}$ we have that
    \begin{equation*}
        f(\w{x}_{k+1}) \leq f(\w{x}_k) - \frac{1}{2} \alpha_k r_k^2.
    \end{equation*}
\end{lemma}

\begin{proof}
    From \cref{ass: lipschitz cont} we have that
    \begin{equation*}
        f(\w{x}_{k+1}) \leq f(\w{x}_{k}) + \langle  \nabla f(\w{x}_{k}), \w{x}_{k+1} - \w{x}_{k} \rangle + \frac{1}{2} \langle  \nabla^2 f(\w{x}_{k}) (\w{x}_{k+1} - \w{x}_{k}), \w{x}_{k+1} - \w{x}_{k} \rangle + \frac{L}{6} r_k^3.
    \end{equation*}
    Using \cref{eq: identity 2} and \cref{eq: identity 3} we take $        f(\w{x}_{k+1}) - f(\w{x}_{k}) $
    \begin{align*}
        & \leq \langle   \w{R}_k \nabla f(\w{x}_{k}), \w{R}_k (\w{x}_{k+1} - \w{x}_{k}) \rangle \\
        & \qquad+ \frac{1}{2}  \langle  \w{R}_k \nabla^2 f(\w{x}_{k})\w{P}_k \w{R}_k(\w{x}_{k+1} - \w{x}_{k}), \w{R}_k (\w{x}_{k+1} - \w{x}_{k}) \rangle  + \frac{L}{6}r_k^3 \\
        & \stackrel{(\ref{eq: identity 1})}{=} - \alpha_k r_k^2
        - \frac{1}{2}\langle (\w{B}_{h,k} + \alpha_k \w{I}_n)\w{R}_k (\w{x}_{k+1} - \w{x}_{k}), \w{R}_k (\w{x}_{k+1} - \w{x}_{k})\rangle  \\
        &  \quad \ \  + \frac{1}{2}  \langle  (\w{R}_k \nabla^2 f(\w{x}_{k})\w{P}_k - \w{B}_{h,k}) \w{R}_k(\w{x}_{k+1} - \w{x}_{k}), \w{R}_k (\w{x}_{k+1} - \w{x}_{k}) \rangle + \frac{L}{6} r_k^3 \\
        & \leq - \alpha_k r_k^2 + \frac{1}{2}  \langle  (\w{R}_k \nabla^2 f(\w{x}_{k})\w{P}_k - \w{B}_{h,k}) \w{R}_k(\w{x}_{k+1} - \w{x}_{k}), \w{R}_k (\w{x}_{k+1} - \w{x}_{k}) \rangle + \frac{L}{6} r_k^3 \\
        & \stackrel{(\ref{eq: definition s_h})}{\leq}
         - \alpha_k r_k^2 + \frac{1}{2}  s_h r_k^2+ \frac{L}{6} r_k^3\\
         & \stackrel{(\ref{ineq: upper bound r_k})}{\leq} 
        - \frac{1}{2} \alpha_k r_k^2 - \frac{1}{2} \left(\alpha_k - s_h - \frac{L\|\w{R}_k \nabla f(\w{x}_{k})\|}{3 \alpha_k} \right) r_k^2 \\
        & \leq - \frac{1}{2} \alpha_k r_k^2.
    \end{align*}
    where the last inequality follows from our assumption on $\alpha_k$. 
\end{proof}

Notice that the reduction in the objective function depends now on $r_k$, in contrast to \cref{lemma: f reduction lambda} which, in absence of \cref{def: P}, uses $\hat{\lambda}(\w{x}_k)$. Therefore, it is more sensible for \cref{alg: inner loop non convex}, under \cref{def: P}, to run the inner loop until the result of \cref{lemma: f reduction definition R} is satisfied.
\red{
\begin{assumption} \label{ass: probability delta}
    There exists $\delta,  \hat{\mu} \in (0,1)$ such that for $\w{R}_k \in \R^{n \times N}$ from Definition 1 it holds $    \| \w{R}_k \nabla f (\w{x}_k) \| \geq \hat{\mu}  \| \nabla f (\w{x}_k) \| \ \text{and}
        \   \| \w{R}_k \nabla f (\w{x}_{k+1}) \| \geq \hat{\mu}  \| \nabla f (\w{x}_{k+1}) \|$ with probability $\delta$.
\end{assumption}
}
\blue{
\cref{ass: probability delta} is typical in subspace methods and trivially satisfied \cite{cartis2022randomised}. In our case, it effectively activates the conditions in \cref{eq: conditions new} during the entire process with some probability $\delta$. In general one should expect a trade-off between the magnitude of $\hat{\mu}$ and $\delta.$ The larger $\hat{\mu}$ becomes, the smaller the probability $\delta$ is expected to be. However, our goal is to apply the algorithm to problems such that $\hat{\mu}$ and $\delta$ are reasonably large (we will see later that this implies fast convergence rates). Of course, we can always achieve this when $n$ approaches $N$, nevertheless this yields large coarse models and hence expensive iterations that we would like to avoid. On the other hand, since $\w{R}$ draws samples from the gradient, it is reasonable to expect large $\hat{\mu}$ and $\delta$ when the gradient is sparse (several zero entries). The following simple example indicates that the proposed algorithmic framework should work well in the case where the Hessian of the problem has low rank. 
}

\blue{
\textbf{Example:} Let $\w{Q} \in \R^{N \times N}$ with $\operatorname{rank} \w{Q} =  r < N$, and $f(\w{x}) \coloneqq \frac{1}{2} \w{x}^T \w{Q} \w{x}$, where $\w{x} = [x_1 \ x_2 \ \cdots \ x_N]^T$. The quadratic function is important in optimization because many of the practical applications are locally quadratic in a neighborhood of a minimizer. Without loss of generality we can assume that $\w{Q}$ is a diagonal matrix with $\{\lambda_i\}_{i = 1}^N$ diagonal elements. Then exactly $r$ elements of $\w{Q}$ are non-zero. Assume without loss of generality they are located at the first $r$ rows. We also have that $\lambda_i = 0$ for all $i = r+1, \ldots, N$. Additionally, we have $\nabla^2 f(\w{x}) = \w{Q}$, i.e., low-rank Hessian matrix, and $ \nabla f(\w{x})  = [\lambda_1 x_1 \ \lambda_2 x_2 \  \cdots \ \lambda_r x_r \ 0 \ \cdots 0]^T$. We wish to compute the probability that, given $\w{R}$ in \cref{def: P}, it holds 
$\| \w{R} \nabla f (\w{x}) \|  =  \| \nabla f (\w{x}) \|$, i.e., $\hat{\mu} = 1$. Assume that the coarse model dimensions are chosen such that $n > r$. Then we wish to compute the probability that $\w{R} \nabla f (\w{x})$ contains the first $r$ entries of the gradient vector. The total number of $n$-samples that contain the first $r$ elements can be found by computing the total number of $n-r$ combinations from the remaining $N-r$ elements (since the first $r$ are known), which is $\binom{N-r}{n-r}$. Then the required probability is
\begin{equation*}
    p = \frac{\binom{N-r}{n-r}}{\binom{N}{n}}.
\end{equation*}
Clearly a large $p$ will be attained if $r$ is small. Then in our example $\hat{\mu}$ and $\delta$ will be large for problems with $r \ll N$. In summary, since \cref{ass: probability delta} requires two successive draws with $\| \w{R} \nabla f (\w{x}) \|  =  \| \nabla f (\w{x}) \|$, we conclude that, \cref{ass: probability delta} is satisfied with probability $\delta = p^2$. Note that this probability is computed for any $\hat{\mu} \in (0,1]$ and if $\hat{\mu} < 1$ then the probability $\delta$ will be larger.
}
\begin{corollary} \label{cor: f reduction non convex probability}
        Let Assumptions \ref{ass: lipschitz cont} and \ref{ass: probability delta} hold and  define $\w{R}_k$ as in \cref{def: P}. Suppose also that the sequence $(\w{x}_k)_{k \in \N}$ is generated by \cref{eq: update rule non convex coarse}. Then for any $\alpha_k \geq s_h + \sqrt{\frac{L \| \w{R}_k \nabla f(\w{x}_k)\|}{2}}$ we have that
    \begin{equation*}
        f(\w{x}_{k+1}) \leq f(\w{x}_k) - \frac{ \hat{\mu}^2 \| \nabla f(\w{x}_{k+1}) \|^2}{8 \alpha_k },
    \end{equation*}
    with probability $\delta$.
\end{corollary}

\begin{proof}
    Squaring the result of \cref{lemma: lower bound rk} and using \cref{ass: probability delta}, we obtain
    \begin{equation*}
        r_k^2 \geq \frac{\| \w{R}_k \nabla f(\w{x}_{k+1}) \|^2}{4 \alpha_k^2} \geq      
        \frac{ \hat{\mu}^2 \| \nabla f(\w{x}_{k+1}) \|^2}{4 \alpha_k^2}.
    \end{equation*}
\end{proof}

We recall the following definitions from \cref{sec: convergence analysis}: $\q{g}_k \coloneqq \| \nabla f (\w{x}_k) \|$ and $\q{g}_k^* \coloneqq \min_{0\leq i \leq k} \q{g}_i$. The next result presents a probabilistic convergence rate when sequence are always generated from (\ref{eq: update rule non convex coarse}) with \cref{def: P}.

\begin{theorem} \label{thm: coarse non-convex probability}
    Let Assumptions \ref{ass: lipschitz cont} and \ref{ass: probability delta} hold and $\w{P}$ as in \cref{def: P}. Suppose also that the sequence $(\w{x}_k)_{k \in \N}$ is generated by \cref{eq: update rule non convex coarse} and 
    \begin{equation*}
        \alpha_k = s_h + \sqrt{\frac{L\| \w{R}_k \nabla f(\w{x}_{k}) \|}{2}},
    \end{equation*}
    where $s_h$ is defined in \cref{eq: definition s_h}. Then for every $\varepsilon > 0$ and 
    \begin{equation*}
    K_\varepsilon \coloneqq \left\lfloor \frac{8(f(\w{x}_0) - f^*)}{\mu^2} \left(\sqrt{\frac{L}{2}} \frac{1}{\varepsilon^{\frac{3}{2}}} + \frac{s_h}{\varepsilon^2}
        \right) + 2\ln{\frac{\q{g}_0}{\varepsilon}} \right\rfloor + 1
\end{equation*}
    we have that $\q{g}_k^* < \varepsilon$ for all $k \geq K_\varepsilon$ with probability $\delta^{K_\varepsilon}$.
\end{theorem}

\begin{proof}
    The proof of the theorem follows \cref{thm: coarse non-convex} but now using \cref{cor: f reduction non convex probability}. Since the reduction shown in \cref{cor: f reduction non convex probability} occurs with fixed probability $\delta$, then the desired convergence rate holds with probability $\delta^{K_\varepsilon}$.
\end{proof}

\cref{thm: coarse non-convex probability} is effectively a probabilistic version of \cref{thm: coarse non-convex} since the iterates are random variables due to $\w{R}_k$. It shows the probability of a single run of \cref{alg: non convex} to achieve the desired tolerance while always performing coarse step. Moreover, we note that we can recover the above convergence rate by setting $\omega=1$ in \cref{thm: coarse non-convex}, even though the analysis is different in this section. As a result, \cref{thm: coarse non-convex} shows a slower convergence of \cref{alg: non convex} when $\omega > 1$.

\subsubsection{The Convex Case}

Assume that \cref{alg: convex} always builds sequences by following the first equation in \eqref{eq: update rule convex} using $\w{R}_k$ from \cref{def: P}. We will provide a probabilistic convergence rate for general convex functions. Similar to the non-convex case it holds
\begin{equation} \label{ineq: upper bound r_k convex}
    r_k \leq \| (\w{R}_k \nabla f(\w{x}_{k}) \w{P}_k + \alpha_k \w{I}_n)^{-1}\| \| \w{R}_k \nabla f(\w{x}_k)\| \leq 
    \frac{\| \w{R}_k \nabla f(\w{x}_k)\|}{\alpha_k},
\end{equation}
for any $\alpha_k > 0$. The next results trivially follow from our analysis in \cref{sec: non-convex probabilities}.
\begin{lemma} \label{lemma: lower bound rk convex}
    Let Assumptions \ref{ass: lipschitz cont} and $\w{R}_k$ as \cref{def: P} hold and suppose that the sequence $(\w{x}_k)_{k \in \N}$ is generated by \cref{eq: update rule non convex coarse}. If $\alpha_k \geq  \sqrt{\frac{L \| \w{R}_k \nabla f(\w{x}_k)\|}{2}}$, then
    \begin{equation*}
        \| \w{R}_k \nabla f(\w{x}_{k+1})\| \leq 2 \alpha_k r_k.
    \end{equation*}
\end{lemma}

\begin{proof}
    Applying the steps in the proof of \cref{lemma: lower bound rk} with $s_h = 0$ and $\w{B}_{h, k} = \w{R} \nabla^2 f_H (\w{x}_{H,k}) \w{P}$, we obtain the desired result.
\end{proof}


\begin{lemma} \label{lemma: f reduction convex definition R}
    Let Assumptions \ref{ass: lipschitz cont} and form $\w{R}_k$ from \cref{def: P}. Suppose also that the sequence $(\w{x}_k)_{k \in \N}$ is generated by the first equation in \cref{eq: update rule convex}. Then for any $\alpha_k \geq \sqrt{\frac{L \| \w{R}_k \nabla f(\w{x}_k)\|}{2}}$ we have that
    \begin{equation*}
        f(\w{x}_{k+1}) \leq f(\w{x}_k) - \frac{2}{3} \alpha_k r_k^2.
    \end{equation*}
\end{lemma}

\begin{proof}
    From \cref{lemma: f reduction definition R} we have that 
    \begin{equation*}
    f(\w{x}_{k+1}) - f(\w{x}_k) \leq - \alpha_k r_k^2 + \frac{1}{2}  s_h r_k^2+ \frac{L}{6} r_k^3
    \end{equation*}
    Since now $\w{R}_k \nabla^2 f(\w{x}_{k}) \w{P}_k$ is available and $s_h = 0$ we have 
    \begin{align*}
        f(\w{x}_{k+1}) - f(\w{x}_k) 
        & \leq - \alpha_k r_k^2 +  \frac{L}{6} r_k^3\\
        & \overset{(\ref{ineq: upper bound r_k convex})}{\leq} - ( \alpha_k -  \frac{L \| \w{R}_k \nabla f(\w{x}_k) \|}{6 \alpha_k}) r_k^2 \leq - \frac{2}{3} \alpha_k r_k^2
    \end{align*}
    where the last inequality follows from our assumption on $\alpha_k$.
\end{proof}

Therefore, in this setting, \cref{alg: inner loop convex} should exit the inner loop according to the bound in  \cref{lemma: f reduction convex definition R}. Next, we show a probabilistic reduction in the value of the objective function in terms of $\|\nabla f(\w{x}_{k+1})\|$.

\begin{corollary} \label{cor: f reduction convex probability}
    Let Assumptions \ref{ass: lipschitz cont} hold and \cref{ass: probability delta} hold and form $\w{R}_k$ from \cref{def: P}.  Suppose also that the sequence $(\w{x}_k)_{k \in \N}$ is generated by the first equation in \cref{eq: update rule convex}. 
    Then for any $\alpha_k \geq \sqrt{\frac{L \| \w{R}_k \nabla f(\w{x}_k)\|}{2}}$ we have that
    \begin{equation*}
        f(\w{x}_{k+1}) - f(\w{x}_k) \leq - \frac{\hat{\mu}^2\|\nabla f(\w{x}_{k+1})\|}{6 \alpha_k}
    \end{equation*}
    with probability $\delta$.
\end{corollary}

\begin{proof}
    As in proof of \cref{cor: f reduction non convex probability}, we combine \cref{lemma: f reduction convex definition R} and \cref{ass: probability delta} to take
    \begin{equation*}
        f(\w{x}_k) - f(\w{x}_{k+1}) \geq \frac{2}{3} \alpha_k r_k^2 \geq \frac{\|\w{R}_k \nabla f(\w{x}_{k+1})\|^2}{6\alpha_k} \geq
        \frac{ \hat{\mu}^2 \|\nabla f(\w{x}_{k+1})\|^2}{6\alpha_k}
    \end{equation*}
    with probability $\delta$.
\end{proof}

\begin{theorem} \label{thm: rate convex coarse probability}
    Let $f$ be a convex function and Assumptions \ref{ass: lipschitz cont} and \ref{ass: probability delta} hold. Suppose also that $R<\infty$, where $R$ is defined in \cref{def: R}. If the sequence $(\w{x}_k)_{k \in \N}$ is generated by \cref{eq: update rule non convex fine} with $\w{R}_k$ as in \cref{def: P} and 
    \begin{equation*}
        \alpha_k = \sqrt{\frac{\omega^3 L \| \w{R}_k \nabla f(\w{x}_k)\|}{2}},
    \end{equation*}
    then there exists $\hat{K}_0 \in \N$ such that convergence rate of the sequence $(f(\w{x}_k))_{k \in N}$ to $f^*$ is given by
    \begin{equation}
        f(\w{x}_k) - f^* \leq \frac{96^2 R^3 L}{2 \mu^4 \left(1 + \frac{3k}{4(3 - \log_{2}\hat{\mu} )} \right)^2}, \label{ineq: rate thm convex probability}
    \end{equation}
    with probability $\delta$ for each $k \geq \hat{K}_0$. Additionally, for every $\varepsilon > 0$ and for
    \begin{equation*}
        K_\varepsilon \coloneqq \left\lceil  \frac{128(3-\log_{2} \hat{\mu}) \sqrt{R^3 L}}{  \hat{\mu}^2 \sqrt{2 \varepsilon}}  \right\rceil + \hat{K}_0,
    \end{equation*}
    we have that $f(\w{x}_k) - f^* < \varepsilon$, for all $k \geq K_\varepsilon$ with probability $\delta^{K_\varepsilon}$.
\end{theorem}

\begin{proof}
    The proof of the theorem follows the steps \cref{thm: rate convex coarse} and thus we state only the major differences. Let $k \in \N$.
    Using \cref{cor: f reduction convex probability} and the fact that $\| \w{R}_k \nabla f(\w{x}_k)\| \leq \| \nabla f(\w{x}_k)\|$ we have that
    \begin{equation} \label{ineq: final f reduction in theorem prob} 
        f(\w{x}_{k}) - f(\w{x}_{k+1} ) \geq \frac{\hat{\mu}^2 \| \nabla f(\w{x}_{k+1})\|^2}{6 \| \w{R}_k \nabla f(\w{x}_k)\|^{\frac{1}{2}} \sqrt{\frac{ L }{2}}} \geq  \frac{\hat{\mu}^2 \q{g}_{k+1}^2}{6 \sqrt{\frac{ L }{2}}\q{g}_{k}^{\frac{1}{2}}},
    \end{equation}
    with probability $\delta$. Recall the index sets $\mathcal{I}_{\infty} \coloneqq \{ m \in \N: \q{g}_{m+1} \geq \frac{1}{4} \q{g}_m \}$ and $\mathcal{I}_{k} \coloneqq \{ m \in \mathcal{I}_\infty: m \leq k \}$ and let $D_k \coloneqq \frac{3k}{2(3 - \log_{2}\hat{\mu})}$. Since $\hat{\mu} \in (0,1)$, we get $0< D_k \leq \frac{k}{2}$, and thus we get $D_{k+1} \geq D_k$ implying $\lim_{k \rightarrow \infty } D_k = \infty$.
    
    \textbf{Case 1:} $|\mathcal{I}_k| \geq D_k$. Then $\mathcal{I}_\infty$ contains infinitely many elements and thus there exists a subsequence  $(m_\ell)_{\ell \in \N}$ such that $ \mathcal{I}_\infty = \{m_0, m_1, \ldots \}$. Following the steps of \cref{thm: rate convex coarse} and setting $C \coloneqq \frac{\hat{\mu}^2}{96 R^{\frac{3}{2}} \sqrt{\frac{L}{2}}}$, we show that 
        \begin{equation*}
        f(\w{x}_{m_{\ell}}) - f(\w{x}_{m_{\ell+1}} )  \geq \frac{\hat{\mu}^2 }{96 \sqrt{\frac{ L }{2}}} \q{g}_{m_{\ell}}^\frac{3}{2} \overset{(\ref{ineq: convexity})}{\geq} C \left(f(\w{x}_{m_{\ell}}) - f^*  \right)^{\frac{3}{2}}
    \end{equation*}
    with probability $\delta$. The remaining of the proof parallels that of \cref{thm: rate convex coarse}. It holds
    \begin{equation*}
        f(\w{x}_k) - f^* \leq \frac{96^2 R^3 L}{2 \mu^4 \left(1 + \frac{3k}{4(3 - \log_{2}\mu )} \right)^2},
    \end{equation*}
     which is true with probability $\delta$ for each $k \in \N$.

    \textbf{Case 2:} $|\mathcal{I}_k| < D_k$. We consider two cases for $\q{g}_{m}$ as previously. \textbf{(a)}  $i \notin \mathcal{I}_\infty$ which implies $\q{g}_{m} < \frac{1}{4}\q{g}_{m-1}$, and \textbf{(b)} $i \in \mathcal{I}_\infty$ which implies  $\q{g}_{m} \leq \frac{2}{\hat{\mu}}\q{g}_{m-1}$ which is true with probability $\delta$.  Then combining these two inequalities with convexity (\ref{ineq: convexity}) we have that
    \begin{equation*}
       f(\w{x}_k) - f^* \leq \frac{R}{2^{\frac{k}{2}}}\q{g}_0.
    \end{equation*}
    Therefore, the analysis of our cases shows that the algorithm will convergence either with $\mathcal{O}(k^{-2})$ or $\mathcal{O}(2^{-k/2})$ and moreover each iteration will achieve such a decrease with probability $\delta$. Since however the latter is a faster rate, there exists $\hat{K}_0 \in \N$ such that \cref{ineq: rate thm convex probability} is true for each $k \geq \hat{K}_0$ and probability $\delta$.

    Furthermore, let $\varepsilon>0$.  To derive the total number of steps, we first solve for $k$ the inequality \cref{ineq: rate thm convex probability}. Therefore, 
    \begin{equation*}
        \left\lceil  \frac{128(3-\log_{2} \hat{\mu}) \sqrt{R^3 L}}{  \hat{\mu}^2 \sqrt{2 \varepsilon}}  \right\rceil
    \end{equation*}
    are the total number of steps of the method if the rate was $\mathcal{O}(k^{-2})$ from the beginning. Since the method may perform the first $\hat{K}_0$ steps from the exponential rate, then $K_\varepsilon$ follows. Note also that the probability $\delta$ is fixed at each iteration, and thus we reach a tolerance $\varepsilon$ after at least $K_\varepsilon$ iterations with an overall probability $\delta^{K_\varepsilon}$.
\end{proof}

\cref{thm: rate convex coarse probability} is a probabilistic version of \cref{thm: rate convex coarse}. However, note that setting $\omega=1$ in the result of \cref{thm: rate convex coarse}, we obtain a slower convergence rate than when using \cref{def: P}. This is true due to the identities of \cref{lemma: identities}, which hold only under \cref{def: P}. In addition, the faster convergence rate of this section suggest that the exit condition of step \ref{step: exit condition} in \cref{alg: inner loop convex} should be replaced by the bound in \cref{lemma: f reduction definition R}.  Further, note that, if $\hat{\mu} = 1$ in \cref{ass: P}, then we obtain the fastest possible rates of convergence in Theorems \ref{thm: coarse non-convex probability} and \ref{thm: rate convex coarse probability}. Therefore, the proposed scheme should work well for problem structures $\hat{\mu}$ is large. Examples of problems with such structures include, but not limited to, problems where the second-order information is concentrated in the first few eigenvalues or the Hessian is typically low-rank \blue{(see discussion followed \cref{ass: probability delta}; the same is also observed in \cite{tsipinakis2024multilevel})}. Finally, note that if $\hat{\mu} = 1$, then the convergence rate in \cref{thm: rate convex coarse probability} becomes exactly the rate proved in \cite{mishchenko2023regularized}, which also indicates that \cref{def: P} constitutes an attractive choice for generating $\w{R}_k$, not only for its efficiency in practice but also when analyzing multilevel or subspace methods.

\section{Numerical Results} \label{sec: experiments}

\begin{figure*}
\centering
\begin{subfigure}{.48\textwidth}
\centering
  \includegraphics[width=.98\linewidth]{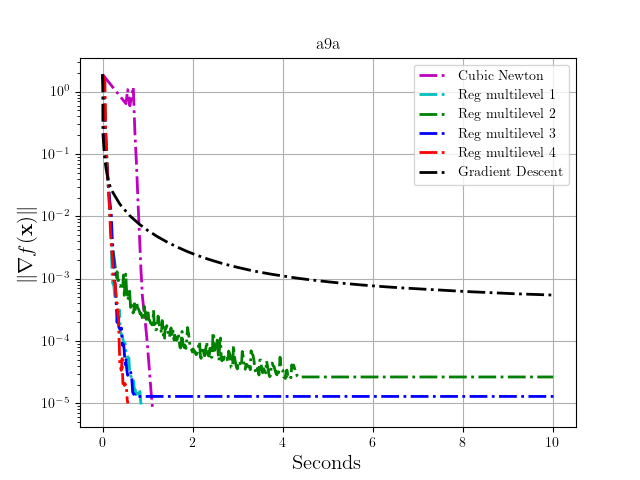}  
\end{subfigure}
\begin{subfigure}{.48\textwidth}
\centering
  \includegraphics[width=.98\linewidth]{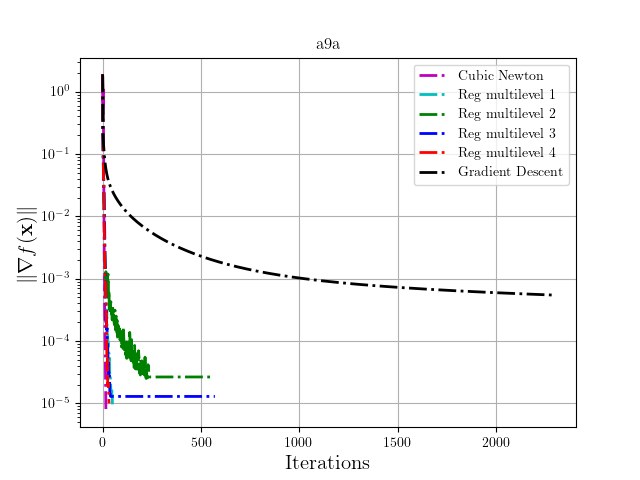}  
\end{subfigure}
\\
\begin{subfigure}{.48\textwidth}
\centering
  \includegraphics[width=.98\linewidth]{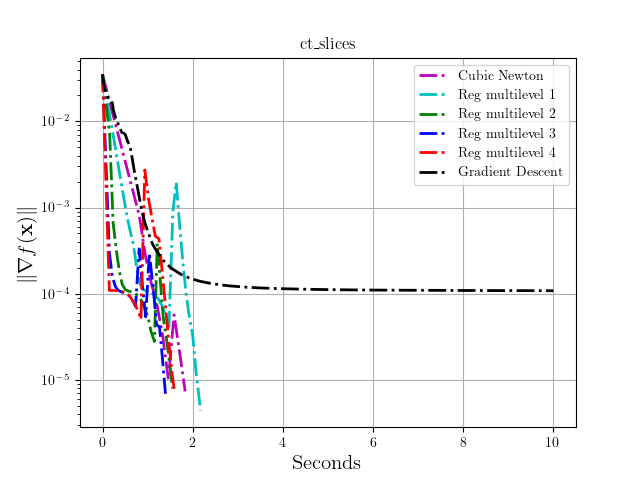}  
\end{subfigure}
\begin{subfigure}{.48\textwidth}
\centering
  \includegraphics[width=.98\linewidth]{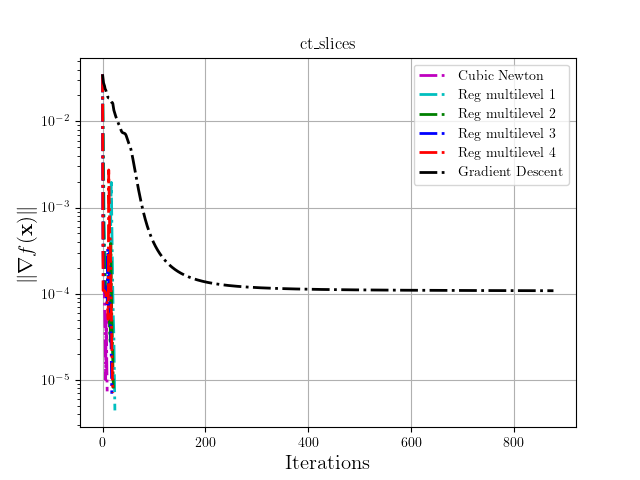}  
\end{subfigure}
\\
\begin{subfigure}{.48\textwidth}
\centering
  \includegraphics[width=.98\linewidth]{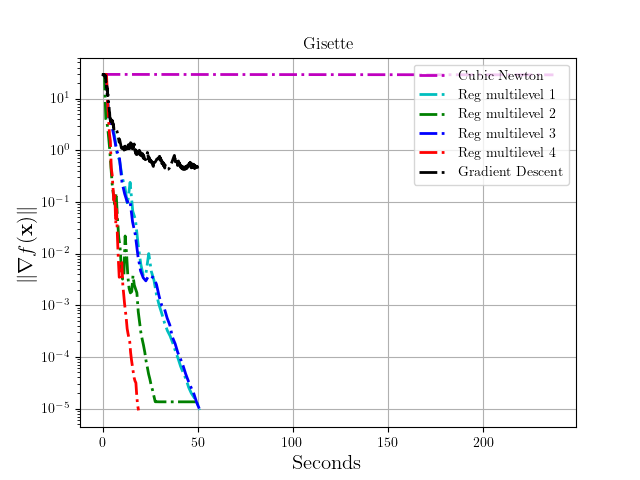}  
\end{subfigure}
\begin{subfigure}{.48\textwidth}
\centering
  \includegraphics[width=.98\linewidth]{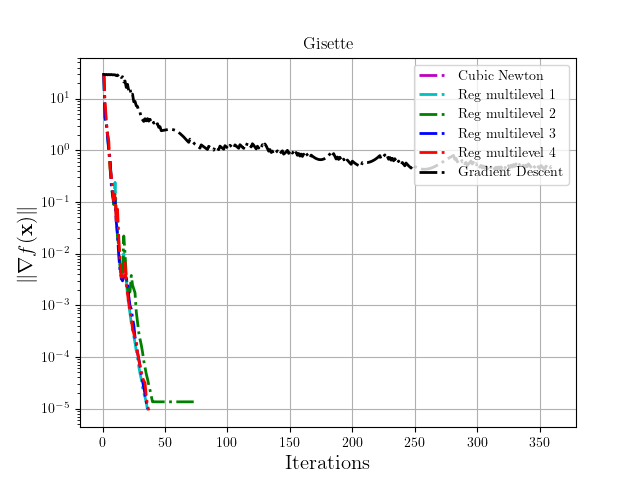}  
\end{subfigure}
\caption{Convergence behaviour of optimization algorithms on logistic regression. The dataset used to produce each plot is mentioned on top.  In all problems, $\w{x}_0$ is initialized randomly in $[0,1]$ from the uniform distribution.
}
\label{fig: logistic regression}
\end{figure*}

In this section we aim to verify our theoretical results through numerical experiments. In particular, we consider real-world problems to demonstrate that the regularized multilevel methods can be significantly faster than Gradient Descent in terms of both number of iterations and wall-clock time. In addition, we plan to show that the rate of regularized multilevel methods can approach the fast rate of Cubic Newton in practice, as it was suggested by Theorems \ref{thm: coarse non-convex} and \ref{thm: rate convex coarse}. We perform comparisons between the regularized multilevel \ref{scenario 1} (i.e., when \eqref{eq:scen1} is used), \ref{scenario 2} (i.e., when \eqref{eq:scen2} is used) and \ref{scenario 3} (i.e., when \eqref{eq:scen3} is used), regularized Newton-type multilevel method (\cref{alg: inner loop convex}), Gradient Descent with Armijo rule for selecting the step-size parameter \cite{nesterov2018lectures} and Cubic Newton \cite{nesterov2006cubic}, on convex and non-convex problems \footnote{All datasets used in the experiments can be found at \url{https://www.csie.ntu.edu.tw/~cjlin/libsvmtools/datasets/}}. 

To obtain the numerical results, we use a python implementation and automatic differentiation \footnote{For a description and efficient implementation of automatic differentiation see remark 2.1 in \cite{tsipinakis2024multilevel}} for computing the gradients and Hessians at each iteration. We will be referring to the different instances of multilevel methods in definitions \ref{scenario 1}, \ref{scenario 2} and \ref{scenario 3} with the name regularized multilevel algorithm $1$, $2$ and $3$, respectively, while we name \cref{alg: inner loop convex} as regularized multilevel algorithm 4. 
In all experiments the coarse dimensions are selected $n = 0.5 N$. Further, we simplify the line search in Algorithms \ref{alg: inner loop convex} and \ref{alg: inner loop non convex} to searching the smallest $i_k \in \N$ such that $\alpha_k \coloneqq 2^{i_k} L_0$ yields the required reductions in the objective function. Additionally, we initialize the line search with $L_0 = 10^{-12}$ for both multilevel and Cubic Newton methods.
\begin{figure*}
\centering
\begin{subfigure}{.47\textwidth}
\centering
  \includegraphics[width=.98\linewidth]{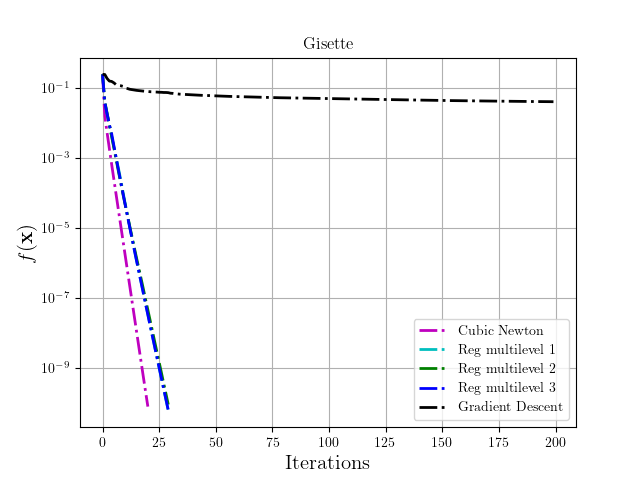}  
\end{subfigure}
\begin{subfigure}{.47\textwidth}
\centering
  \includegraphics[width=.98\linewidth]{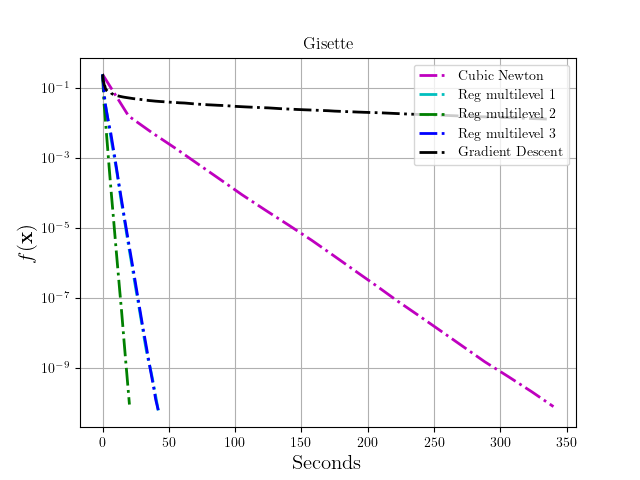}  
\end{subfigure}
\\
\begin{subfigure}{.47\textwidth}
\centering
  \includegraphics[width=.98\linewidth]{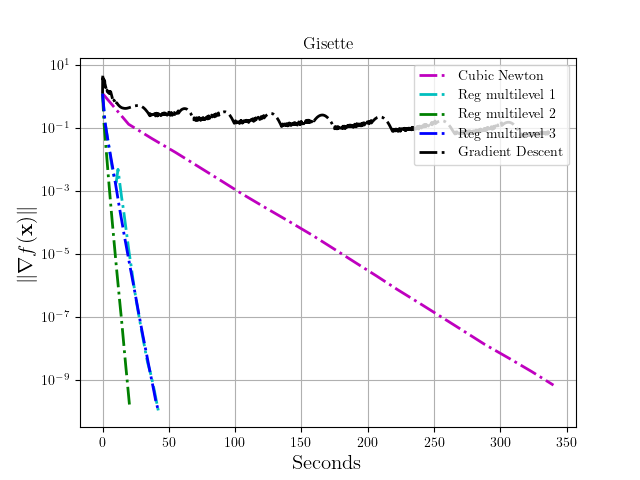}  
\end{subfigure}
\caption{Convergence behaviour of optimization algorithms on non-linear least-squares. We set $\w{x}_0=0$. The regularized multilevel method 2 uses $r=500$.
}
\label{fig: nlls gisette}
\end{figure*}
Given a collection of points of the form $\{\w{a}_i, b_i \}_{i=1 \ldots m}$, first, we are interested in solving the logistic regression problem which has the following form
\begin{equation*}
\min_{\w{x} \in \R^N} \frac{1}{m} \sum_{i=1}^m \log(1 + \exp(- \langle \w{a}_i, \w{x} \rangle)) + \frac{\lambda}{2} \| \w{x}\|^2,   
\end{equation*}
where $\lambda > 0$ (\blue{for this experiment we use $\lambda = 10^{-6}$}). The above problem is convex and thus both algorithms (Alg. \ref{alg: inner loop convex} and \ref{alg: inner loop non convex}) can be applied.
The results appear in \cref{fig: logistic regression}. Clearly, in all experiments, Gradient Descent is significantly slower than all other algorithms, even though it has by far the lowest cost of forming the iterates. 
Moreover, multilevel methods enjoy a similar rate convergence to Cubic Newton, which makes them very efficient in large scale optimization as they generate their iterations in a coarse level. The efficiency of multilevel methods compared to the Cubic Newton is more evident when the dimensions of the original model become large, see e.g. \cref{fig: logistic regression} for the Gisette dataset. In this problem $N=5000$ which means that by the time the Cubic Newton  performs a single iteration the multilevel methods have already converged to a solution with a high accuracy. Among the multilevel methods, the method 4 (generating the iterates by using the first equation in \ref{eq: update rule convex}) is typically the fastest and most consistent to reach the required tolerance. This is expected since the problem is convex. 
\begin{figure*}
\centering
\begin{subfigure}{.48\textwidth}
\centering
  \includegraphics[width=.98\linewidth]{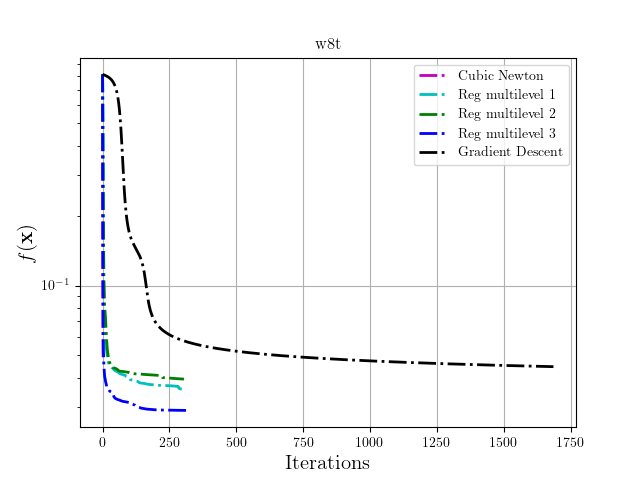}  
\end{subfigure}
\begin{subfigure}{.48\textwidth}
\centering
  \includegraphics[width=.98\linewidth]{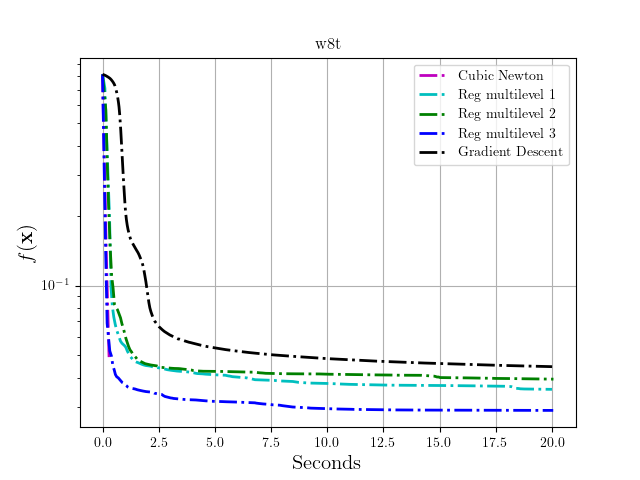}  
\end{subfigure}
\\
\begin{subfigure}{.48\textwidth}
\centering
  \includegraphics[width=.98\linewidth]{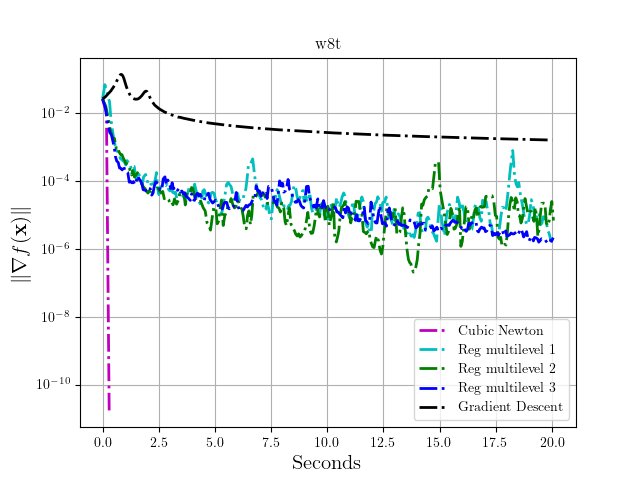}  
\end{subfigure}
\caption{Convergence behaviour of optimization algorithms on non-linear least-squares problem. We set $\w{x}_0=0$. The regularized multilevel method 2 uses $r=60$.
}
\label{fig: nlls w8}
\end{figure*}

Next we wish to solve the non-linear least-squares problems which has the form
\begin{equation*}
\min_{\w{x} \in \R^N} \frac{1}{m} \sum_{i=1}^m (b_i - \phi(\langle \w{a}_i, \w{x} \rangle))^2, \quad \phi(x) \coloneqq \frac{\exp(x)}{1 + \exp(x)}.
\end{equation*}
This problem is non-convex and thus only \cref{alg: inner loop non convex} can be applied. As before, we use  \eqref{eq:scen1},  \eqref{eq:scen2} and \eqref{eq:scen3} for constructing  $\w{B}_{h,k}$. The results appear in Figures \ref{fig: nlls gisette} and \ref{fig: nlls w8}. Figure \ref{fig: nlls gisette} shows the convergence behavior of the algorithms using the Gisette dataset. Reaching the global minimum at this problem is challenging as it has several flat areas and saddle points along the optimization path which pose difficulties to optimization methods that do not utilize enough second-order information. The figure shows that the multilevel methods clearly outperform Cubic Newton and Gradient Descent methods. In particular, Gradient Descent gets stuck at a saddle far from the minimum, while Cubic Newton does not convergence before the time limit is exceeded. Between the multilevel methods, the regularized multilevel method 2 is the fastest method due to the cheaper cost of using a low rank decomposition. On the other hand, regularized multilevel methods 1 and 3 enjoy a faster rate (their performance coincide). Similarly, in Figure \ref{fig: nlls w8}, we observe that Gradient Descent converges to a point that yields higher value in the objective function compared to multilevel methods. On the other hand, Cubic Newton converges rapidly to a local minimum but this point has a higher value than the one computed with multilevel methods. We suspect the lower value of multilevel methods occurs due to the randomness induced by \cref{def: P}. In any case, the multilevel methods are efficient and enjoy a satisfactory performance in this experiment.

\section{Conclusions}

We develop new multilevel methods that are based on the regularization of the Hessian matrix or an appropriate preconditioning matrix. We analyze these methods for arbitrary functions with Lipschitz continuous Hessian matrices and for general convex functions. For the former class of functions we prove a rate of convergence that interpolates between the rate of Gradient Descent and that of the Cubic Newton method, and we describe instances of the method and discuss problem structures such that our rates are expected to be faster than Gradient Descent. For general convex functions we show the fast $\p{O}(k^{-2})$ rate. Our numerical experiments verify the fast rates in real-world problems in both non-convex and convex settings.

\section*{Acknowledgments}
We gratefully thank the anonymous reviewers for their constructive feedback and valuable insights, which have significantly enhanced the quality of this manuscript.
 We also extend our sincere thanks to the editors for their careful handling of the review process. Additionally, the second author acknowledges support from EPSRC Grant (EP/W003317/1), Google Cloud Research Credit Award, and an NVIDIA Academic Grant Award.

\bibliographystyle{siamplain}
\bibliography{references}
\end{document}